\documentclass[11pt,reqno]{amsart}

\usepackage[utf8]{inputenc}

\usepackage{amsmath,amsthm,amssymb,amscd,mathrsfs,amsbsy}
\usepackage{mathtools,color,esint,bm,float}
\usepackage{booktabs}
\usepackage[perpage]{footmisc}
\usepackage[shortlabels]{enumitem}
\usepackage{aliascnt}
\usepackage{hyperref}
\usepackage[capitalize,nameinlink]{cleveref}   
\usepackage{comment}
\usepackage{graphicx}
\usepackage{caption,subcaption}

\definecolor{darkgreen}{rgb}{0.0, 0.6, 0.13}

\newcommand{\Rbb}{\mathbb{R}}

\newcommand{\Nbb}{\mathbb{N}}

\newcommand{\Ccal}{\mathcal{C}}

\newcommand{\Var}{\mathrm{Var}}
\newcommand{\Ubv}{\mathcal{U}_{\mathrm{BV}}}
\newcommand{\UbvN}[1][N]{\mathcal{U}_{\mathrm{BV},#1}}
\newcommand{\Vxi}{\mathrm{V}_{\xi}}
\newcommand{\xtar}{x_{\mathrm{target}}}

\theoremstyle{plain}
\newtheorem{theorem}{Theorem}[section]
\newtheorem*{mainresult}{Main Result}
\newaliascnt{lemma}{theorem}
\newtheorem{lemma}[lemma]{Lemma}
\aliascntresetthe{lemma}
\crefname{lemma}{Lemma}{Lemmas}
\newaliascnt{proposition}{theorem}

\aliascntresetthe{proposition}
\crefname{proposition}{Proposition}{Propositions}
\newaliascnt{corollary}{theorem}
\newtheorem{corollary}[corollary]{Corollary}
\aliascntresetthe{corollary}
\crefname{corollary}{Corollary}{Corollaries}
\theoremstyle{remark}
\newaliascnt{remark}{theorem}
\newtheorem{remark}[remark]{Remark}
\aliascntresetthe{remark}
\crefname{remark}{Remark}{Remarks}
\theoremstyle{definition}
\newaliascnt{definition}{theorem}
\newtheorem{definition}[definition]{Definition}
\aliascntresetthe{definition}
\crefname{definition}{Definition}{Definitions}
\newaliascnt{assumption}{theorem}

\aliascntresetthe{assumption}
\crefname{assumption}{Assumption}{Assumptions}

\title[Optimal Graphon Control]{Optimal Graphon Control of Collective Dynamics}

\author{Joseph K. Miller}
\address{Department of Mathematics, Massachusetts Institute of Technology}
\email{jkm314@mit.edu}

\author{Brandon Su}
\address{Institute for Computational and Mathematical Engineering, Stanford University}
\email{jcbsu@stanford.edu}
\date{}

\begin{document}
\maketitle

\begin{abstract}
We study optimal control of non-linear graphon dynamics in which the interaction kernel itself is the control. Motivated by opinion formation in social networks, we ask how changing who influences whom can steer a population toward a desired outcome. We develop a framework using interaction kernels with bounded variation in the receiver label, providing compactness while accommodating finite networks. We prove existence of an optimal pair and show that finite-agent approximations recover the continuum optimal value under successive refinement of the network and initial data, with subsequential convergence of suitably chosen discrete optimizers. Numerical experiments use neural networks to approximate optimal controls for bounded-confidence and graphon Kuramoto dynamics, illustrating how changes in interactions can prevent fragmentation and promote synchronization.
\end{abstract}

\tableofcontents

\section{Introduction}

Graphon dynamics describe the evolution of a population interacting through a weighted network. A graphon is a continuum analogue of the network's adjacency matrix. With agents indexed by labels in $I=[0,1]$, the equation takes the form
\[
    \partial_t x(t,\xi)
    = \int_I w(t,\xi,\zeta)\,\phi\!\left(x(t,\zeta)-x(t,\xi)\right)\,d\zeta,
    \qquad x(0,\xi)=x_0(\xi).
\]
Here $t\in[0,T]$ is time, $T>0$ is a fixed horizon, and $x(t,\xi)\in\Rbb$ is the state of the agent with label $\xi$, starting from the initial profile $x_0$. The integration variable $\zeta$ labels the agents influencing $\xi$, and $d\zeta$ is the normalized population measure. The kernel $w(t,\xi,\zeta)$ specifies the influence of sender $\zeta$ on receiver $\xi$, while $\phi:\Rbb\to\Rbb$ describes the response to their state difference. We use the term graphon for a bounded, nonnegative interaction kernel, allowing directed influence through kernels that need not be symmetric.

\subsection{Continuum limits and kinetic theory}

This equation arises by replacing sums over increasingly large networks with integrals over a continuum of labels. Medvedev \cite{doi:10.1137/130943741} rigorously justified this passage for dynamical systems on convergent sequences of dense simple and weighted graphs. Ayi and Pouradier Duteil \cite{Ayi2021} subsequently studied graph and mean-field limits for collective dynamics with evolving influence weights. These results provide a way to retain network heterogeneity in the large-population limit; see also the survey \cite{Ayi26}.

The graphon description is closely connected to kinetic theory. For homogeneous interactions, $w\equiv1$, the distribution of states evolves by a Vlasov-type transport equation: agents move through state space under the interaction field generated by the population \cite{Dobrushin}. For a general network, the distribution of states alone usually loses essential information about who influences whom. Keeping the labels produces a family of state distributions, each transported by a field weighted by the graphon. The equation above describes the characteristics when each label carries a single state, or equivalently a point mass in state space. Jabin, Poyato, and Soler \cite{JabinPoyatoSoler} place this connection in a broader mean-field theory for non-exchangeable systems, including sparse and asymmetric networks.

Different choices of $\phi$ describe different collective phenomena. In opinion dynamics, bounded confidence expresses the tendency to respond only to sufficiently similar opinions, as in the Hegselmann--Krause model \cite{Hegselmann-Krause}. In the appendix, we study a smoothed continuous-time variant. We also study graphon Kuramoto dynamics for identical natural frequencies in a rotating frame \cite{Kuramoto1975}. Here $\phi(z)=\sin z$, and the states represent oscillator phases. Related alignment mechanisms occur in the second-order Cucker--Smale flocking model \cite{CuckerSmale,CuckerSmaleProof}, which couples positions and velocities.

\subsection{The control problem}

These models lead to the question: \emph{what happens if the graphon itself is controlled?} In social networks, a platform can influence opinion formation by changing who hears whom through recommendations and feeds. For a given initial state, the controller chooses $w(t,\xi,\zeta)$ as a measurable function of time and the sender and receiver labels, without feedback dependence on the evolving states. Each receiver has a fixed attention budget, expressed by $\int_I w(t,\xi,\zeta)\,d\zeta=1$, and individual influence is bounded by a fixed cap.

The controller minimizes the cost
\[
    J[x,w] := \int_0^T \Bigl(\Var(x(t)) + \beta\|w(t)\|_{L^2(I^2)}^2\Bigr)\,dt
    + \gamma\bigl|\overline{x}(T)-x_{\text{target}}\bigr|^2.
\]
Here $\overline{x}(t):=\int_I x(t,\xi)\,d\xi$ and $\Var(x(t)):=\int_I|x(t,\xi)-\overline{x}(t)|^2\,d\xi$ denote the population mean and variance. The weights $\beta,\gamma\ge0$ and target $x_{\text{target}}\in\Rbb$ are prescribed. The three terms penalize disagreement over time, concentration of attention, and deviation of the final mean from the target, respectively.

\subsection{Related work and contribution}

Related graphon control problems differ in the role assigned to the interaction network. Gao and Caines \cite{GaoCaines2020} study external control of linear network systems with prescribed graphon operators, including minimum-energy and linear-quadratic problems and their finite-network approximations. Cao and Lauri\`ere \cite{CaoLauriere2025} develop a probabilistic theory of graphon mean field control in which agents' dynamics are controlled while the interaction graphon is \emph{prescribed}. Graphon games likewise address strategic optimization in heterogeneous populations with a fixed interaction network \cite{doi:10.1137/24M1709339}. More directly, Djete \cite{Djete2025} studies controlled interactions in non-exchangeable stochastic systems, using a relaxed-control framework under a nondegeneracy assumption on the diffusion.

To solve the deterministic nonlinear control problem above, we develop a framework based on bounded variation (BV) of the interaction kernel with respect to the receiver label. This regularity provides compactness for existence while accommodating finite-agent approximations. Our main result, stated precisely in \Cref{thm:main}, is the following:
\begin{mainresult}
Assume that $\phi$ is bounded and globally Lipschitz and that $x_0\in L^2(I)$. For nonnegative, row-stochastic, uniformly bounded controls satisfying a fixed receiver-variation bound, the following hold:
\begin{enumerate}
    \renewcommand{\labelenumi}{(\roman{enumi})}
    \item The continuum problem admits an optimal pair $(x^\star,w^\star)$.
    \item For piecewise-constant initial data, the finite-agent optimal values converge as the dyadic network is refined. Subsequently refining the initial data recovers the continuum optimum for $x_0$.
    \item Along suitable dyadic subsequences, finite-agent optimal controls converge weakly in $L^2$, and their states converge strongly in $L^2$, to a continuum optimal pair.
\end{enumerate}
\end{mainresult}

\Cref{sec:setting} sets out the control problem and its finite-agent discretization. \Cref{sec:existence} proves existence for finite networks, and Section~4 proves the main result. The numerical appendix, \Cref{sec:numerics}, follows the neural-network approach of Carmona and Lauri\`ere \cite{carmona2023deep} to learn candidate controls by approximately minimizing the discretized cost for bounded-confidence and Kuramoto dynamics. These experiments illustrate how changing interactions can prevent fragmentation or promote synchronization; receiver variation is measured after training, rather than constrained to a prescribed bound during optimization.

\subsection{Acknowledgements} 
The authors would like to thank Elina Chatzidimitriou for pointing out the paper \cite{carmona2023deep} during her candidacy examination. Additionally, the first author would like to thank Natalie Ayi for discussions on graphon dynamics which got this project started.

The authors used AI tools developed by Anthropic and OpenAI for mathematical discussions and assistance with editing this manuscript. They take full responsibility for the mathematical arguments, references, and final text.

\section{Setting and Notations}\label{sec:setting}

\subsection{The Graph-Limit Dynamics}

Let \(I := [0,1]\) be a set of labels. We consider the \emph{graph-limit} (or continuum) dynamics (see \cite{Ayi26} for a recent survey), where each agent is indexed by a label \(\xi \in I\), and the state of the agent at time \(t\) is denoted by $x(t,\xi) \in \Rbb$.

The deterministic dynamics are given by
\begin{equation}\label{eqn:graphon-dynamics}
\partial_t x(t,\xi) = \int_I w(t,\xi,\zeta)\,
\phi\!\left(x(t,\zeta)-x(t,\xi)\right) \, d\zeta, \qquad x\in \mathcal{C}^1\Big([0,T] \rightarrow L^2(I; \Rbb)\Big)
\end{equation}
where \(w(t,\xi,\zeta)\) is a possibly time-dependent graphon encoding interaction strengths and \(\phi : \Rbb \to \Rbb\) is the interaction function.

\subsection{Graphon Controlled Dynamics}

We now introduce \emph{control} by allowing the graphon to be time-dependent and chosen by a controller. We restrict to the case where the controller can depend on time and on the sender and receiver labels, but \emph{not} on the current state. We consider the controlled dynamics \eqref{eqn:graphon-dynamics}, where the graphon is now the control: $w : [0,T] \rightarrow \mathcal{U}$ with
\begin{equation}\label{eqn:asymmetric-graphon}
    \mathcal{U} = \left\{ w\in L^\infty(I^2) \, \Big| \, \int_{I} w(\xi,\zeta) \, d\zeta = 1, \ 0 \leq w(\xi,\zeta) \leq \Lambda \text{ for a.e. } \xi,\zeta \in I\right\}.
\end{equation}

For each receiver $\xi \in I$, the row $w(\xi,\cdot)$ is a probability density over senders $\zeta \in I$, so $w$ is a row-stochastic attention kernel. The constant $\Lambda \ge 2$ is a fixed \emph{attention-saturation cap}: it forbids any single sender from absorbing more than $\Lambda \, d\zeta$ of the receiver's attention. Any fixed $\Lambda > 1$ works for all sufficiently fine discretization levels, and we take $\Lambda \ge 2$ so that the discrete classes $\mathcal{U}_N$ below are non-empty for every $N \ge 2$. We define $t \mapsto w(t)$ as \emph{admissible} if it is measurable and $w(t) \in \mathcal{U}$ for a.e. $t \in [0,T]$.
\begin{definition}\label{def:U-BV}
For $w \in L^\infty(I^2)$, the \emph{receiver-label total variation} is the essential variation of the curve $\xi \mapsto w(\xi,\cdot) \in L^2(I)$,
\begin{equation}\label{eqn:receiver-tv}
\Vxi(w)
:= \sup\left\{\,
\sum_{k=1}^{m}
\left(\int_I \bigl|w(\xi_k,\zeta) - w(\xi_{k-1},\zeta)\bigr|^2 \, d\zeta \right)^{1/2}
\;:\; 0 = \xi_0 < \xi_1 < \cdots < \xi_m = 1
\,\right\},
\end{equation}
where the supremum is taken over the precise representative of $w$. For a fixed radius $R> 0$, set
\begin{equation}\label{eqn:U-BV-def}
  \Ubv
  := \bigl\{\, w \in \mathcal{U} \;:\; \Vxi(w) \le R \,\bigr\},
\end{equation}
where $\mathcal{U}$ is defined in \eqref{eqn:asymmetric-graphon}. A control $t \mapsto w(t)$ is \emph{$\mathrm{BV}$-admissible} if it is measurable and $w(t) \in \Ubv$ for a.e.\ $t \in [0,T]$.
\end{definition}

Equivalently using the theory of bounded variation metric space value functions \cite[Definition 2.1, Remark 2.2]{Ambrosio1990metricBV}, $\Vxi(w) = |D_\xi w|(I)$ is the total mass of the $L^2(I_\zeta)$-valued
measure $D_\xi w$, and $\Vxi(w) < \infty$ if and only if
$\xi \mapsto w(\xi,\cdot) \in \mathrm{BV}\bigl(I; L^2(I)\bigr)$ .

\begin{remark}\label{rmk:bv-not-lipschitz}
For our purposes, the class of controls must be $\mathrm{BV}$ in $\xi$, \emph{not} Lipschitz in $\xi$. Under a uniform $W^{1,\infty}_\xi$ constraint\footnote{Here $W^{1,\infty}_\xi:=W^{1,\infty}(I_\xi;L^2(I_\zeta))$ is the Sobolev space of $L^2(I_\zeta)$-valued functions.}, the recovery argument of \Cref{lem:liminf} fails, because $P_N w$ is piecewise constant in $\xi$ and hence not Lipschitz, so no recovery sequence exists inside the class. $\mathrm{BV}$ is the natural choice because $P_N$ is variation-diminishing and therefore preserves it. 
\end{remark}

The goal of the control is to steer the system toward \emph{consensus} or \emph{alignment}. We track the running variance and the mean by
$$
\Var(x) := \int_I \Vert x(\xi) - \overline{x} \Vert^2 \, d\xi, \qquad
\overline{x} := \int_I x(\xi) \, d\xi.
$$
Fix weights $\beta, \gamma \geq 0$ and a target $x_{\text{target}} \in \Rbb$. For an admissible control $w \in \mathcal{U}$, define the cost
\begin{equation}\label{eqn:J-continuum}
J[x, w] := \int_0^T \big(\Var(x(t)) + \beta \Vert w(t) \Vert^2_{L^2(I^2)}\big) \, dt + \gamma \vert \overline{x}(T) - x_{\text{target}} \vert^2. 
\end{equation}
We can alternatively write the cost as directly depending on the initial data $x_0 \in L^2(I)$ and admissible control $w: [0,T] \rightarrow \mathcal{U}$ by defining
\begin{equation}\label{eqn:J-continuum-alt}
    \mathcal{J}(w; x_0) := J[S(w; x_0), w], 
\end{equation}
where $S(w;x_0)$ is the solution map of \eqref{eqn:graphon-dynamics}. This naturally leads to the optimal control problem given by
\begin{equation}\label{eqn:A-infty-y-def}
    A_\infty(x_0) := \inf \big\{ \mathcal{J}(w; x_0) \, : \, w \in L^\infty([0,T]; \Ubv) \big\}.
\end{equation}

\subsection{The Finite Agent System}
In this section, we define the $N$-agent dynamics. Let $\{\xi_i\}_{i=1}^N \subset I$ be a discretization of labels given by $\xi_i = (i - 0.5)/N$.  Then, we approximate the controlled dynamics for a time-and-label control $w = w(t, \xi, \zeta)$ using the histogram function approximation 
\begin{equation}\label{eqn: finite agent approximation}
    x_N(t, \xi) = \sum_{i=1}^N x_N^i(t) \, \mathbf{1}_{I_i^N}(\xi), \qquad I_i^N = \left[\frac{i-1}{N}, \frac{i}{N}\right),
\end{equation}
where $\{x_N^i(t)\}_{i=1}^N \subset \Rbb$ is the list of the states of $N$ agents. We discretize the graphon slice-wise in time by its cell averages
\begin{equation}\label{eqn: finite graphon approximation B}
    \overline{w}_N^{i,j}(t) = N^2 \int_{I_i^N\times I_j^N} w(t, \xi,\zeta)\, d\xi \, d\zeta,
\end{equation}
so that the piecewise-constant kernel $w^N(t) := \sum_{i,j=1}^N \overline{w}_N^{i,j}(t) \,\mathbf{1}_{I_i^N} \otimes \mathbf{1}_{I_j^N}$ is exactly $P_N w(t)$ for the projection $P_N$ of \Cref{def:PN} below; on the discretized labels we write $\overline{w}_N(t, \xi_i, \zeta_j) = \overline{w}_N^{i,j}(t)$.

\begin{definition}[Projection operators]\label{def:PN}
For dyadic $N$, define the projection operators $P_N$ acting on states $x \in L^2(I)$ and on kernels $w \in L^\infty(I^2)$ by
\begin{equation}\label{eqn:PN-state}
    P_N x(\xi) = \sum_{i=1}^N x_N^i \, \mathbf{1}_{I_i^N}(\xi),
    \qquad
    x_N^i = N \int_{I_i^N} x(\zeta)\, d\zeta,
\end{equation}
\begin{equation}\label{eqn:PN-kernel}
    P_N w(\xi,\zeta) = \sum_{i,j=1}^N \overline{w}_N^{i,j} \,
    \mathbf{1}_{I_i^N}(\xi) \mathbf{1}_{I_j^N}(\zeta),
    \qquad
    \overline{w}_N^{i,j} = N^2 \int_{I_i^N \times I_j^N} w(\xi,\zeta)\, d\xi\, d\zeta.
\end{equation}
Note that $P_N P_M = P_M$ whenever $M \le N$ and both are dyadic.
\end{definition}

Throughout, $M$ and $N$ denote dyadic integers, $M = 2^m$ and $N = 2^n$ with $m \le n$, and we write $x_N := P_N x$ and $w_N := P_N w$ for the projections of \Cref{def:PN}. We can alternatively stack $N$ agents and form the vector
$$
\vec{x}_N(t) = (x_N^1(t), x_N^2(t), \cdots, x_N^N(t)) \in \Rbb^N
$$
and the control $w^N(t, \xi_i, \zeta_j)$ takes value from the admissible set
\begin{equation}\label{eqn: U nonsymmetric}
\mathcal{U}_N = \Big\{w \in [0,\Lambda]^{N \times N}: \sum_{j=1}^{N} w_{ij} = N, \text{ for all } i = 1, \ldots, N \Big\}.
\end{equation}
\begin{remark}
The row-sum normalization $\sum_j w_{ij} = N$ is the discrete counterpart of the continuum row-stochasticity $\int_I w(\xi,\zeta)\,d\zeta = 1$: the block-average embedding $\overline{w}_N^{i,j} = N^2 \int_{I_i^N \times I_j^N} w$ of \eqref{eqn: finite graphon approximation B} maps $\mathcal{U}$ into $\mathcal{U}_N$, and conversely the histogram kernel of a matrix in $\mathcal{U}_N$ satisfies $0 \le w_N \le \Lambda$ and $\int_I w_N(\xi,\zeta)\, d\zeta = 1$. We parametrize the row-stochastic form $W_{ij} := w_{ij}/N$ by a row-softmax normalized neural network in \Cref{sec:nn-arch}.
\end{remark}

\begin{remark}\label{rmk:no-self-loop}
    We deliberately do \emph{not} impose the zero-diagonal constraint $w_{ii} = 0$ familiar from finite interacting-agent models. The reason is that vanishing on the  diagonal is incompatible with a uniform bound on the receiver variation $\Vxi$ of \eqref{eqn:receiver-tv-discrete}. For instance, the uniform zero-diagonal matrix $w_{ij} = \frac{N}{N-1}(1 - \delta_{ij})$ has $\Vxi = \sqrt{2N}\,(1 + O(1/N))$. Imposing $w_{ii} = 0$ would therefore obstruct the approximation of continuum controls by discrete ones. However, allowing self-loops is harmless: the self-interaction enters the dynamics only through the term $\frac{1}{N} w_{ii}\, \phi(0)$.
\end{remark}
At level $N$, we also define 
\begin{equation}\label{eqn:U-BV-N-def}
  \UbvN
  :=  \{ w\in \mathcal{U}_N : \Vxi(w) \leq R \}
\end{equation}
where we are defining the variation as the discrete sum
\begin{equation}\label{eqn:receiver-tv-discrete}
  \Vxi(w)
  = \sum_{i=1}^{N-1}
    \left( \frac{1}{N} \sum_{j=1}^{N} (w_{i+1,j} - w_{i,j})^2 \right)^{1/2}.
\end{equation}
Fix a target state $\xtar \in \Rbb$. We can write the $N$-agent cost functional as a Bolza-problem  type \cite{Cesari1983} by defining for measurable $w : [0,T] \rightarrow \mathcal{U}_N$ and $\vec{x}_N : [0,T] \rightarrow \Rbb^N$,
\begin{equation}\label{eqn: JN definition}
J_N[\vec{x}_N, w] =  \int_0^T \Big[\Var(\vec{x}_N(s)) + \beta \frac{1}{N^2} \sum_{i,j=1}^N |w_{i,j}(s)|^2 \Big] \, ds + \gamma | \overline{\vec{x}}_N(T) - \xtar|^2.
\end{equation}
In the above equation, $\Var(\vec{x}_N(s))$ is the standard variation of a list of numbers.

The next lemma collects the two properties of $\Vxi$ used repeatedly below: the continuum variation of \Cref{def:U-BV} agrees with the discrete sum \eqref{eqn:receiver-tv-discrete} on histogram kernels, and $\Vxi$ controls (indeed, is characterized by) translations in the receiver label.

\begin{lemma}\label{lem:variation-properties}
Let $\Vxi$ be as in \Cref{def:U-BV}.
\begin{enumerate}
    \item[(i)] \emph{(Histogram consistency.)} Let $w = \sum_{i,j=1}^N \overline{w}_N^{i,j} \mathbf{1}_{I_i^N}(\xi) \mathbf{1}_{I_j^N}(\zeta)$ be a piecewise-constant kernel. Then
    \begin{equation}\label{eqn:V-xi-histogram}
        \Vxi(w) = \sum_{i=1}^{N-1} \left( \frac{1}{N} \sum_{j=1}^{N} (\overline{w}_N^{i+1,j} - \overline{w}_N^{i,j})^2 \right)^{1/2},
    \end{equation}
    i.e.\ the continuum essential variation coincides with the discrete sum \eqref{eqn:receiver-tv-discrete}.
    \item[(ii)] \emph{(Translation estimate and characterization.)} For $w \in L^\infty(I^2)$, viewed as the curve $\xi \mapsto w(\xi,\cdot) \in L^2(I_\zeta)$, and any $h \in (0,1)$,
    \begin{equation}\label{eqn:V-xi-translation}
        \int_0^{1-h} \bigl\| w(\xi + h, \cdot) - w(\xi,\cdot) \bigr\|_{L^2(I_\zeta)} \, d\xi \;\le\; h \, \Vxi(w),
    \end{equation}
    and moreover
    \begin{equation}\label{eqn:V-xi-sup-char}
        \Vxi(w) = \sup_{0 < h < 1} \frac{1}{h} \int_0^{1-h} \bigl\| w(\xi + h, \cdot) - w(\xi,\cdot) \bigr\|_{L^2(I_\zeta)} \, d\xi .
    \end{equation}
\end{enumerate}
\end{lemma}

\begin{proof}
The argument relies on the classical $\mathrm{BV}$ increment estimate of Ambrosio \cite{Ambrosio1990metricBV} and the lower semicontinuity result of Ambrosio, Fusco, and Pallara \cite{AmbrosioFuscoPallara2000}, with the Hilbert-valued extension explained below. Set $H=L^2(I_\zeta)$ and $g(\xi)=w(\xi,\cdot)$, and let $V_H(g;K)$ denote the essential variation of $g$ on an interval $K$.

\emph{(i)} For a histogram kernel, define
\[
    v_i:=\sum_{j=1}^N\overline{w}_N^{i,j}\mathbf{1}_{I_j^N}(\cdot)\in H,
    \qquad i=1,\ldots,N.
\]
Then $g(\xi)=v_i$ on each cell $I_i^N$, and its distributional derivative consists of the jumps at the interior cell boundaries:
\[
    Dg=\sum_{i=1}^{N-1}(v_{i+1}-v_i)\delta_{i/N}.
\]
Since these atoms have distinct locations and $|I_j^N|=1/N$, their total variation is
\[
    \Vxi(w)=\sum_{i=1}^{N-1}\|v_{i+1}-v_i\|_H
    =\sum_{i=1}^{N-1}\left(\frac{1}{N}\sum_{j=1}^N
    |\overline{w}_N^{i+1,j}-\overline{w}_N^{i,j}|^2\right)^{1/2},
\]
which proves \eqref{eqn:V-xi-histogram}.

\emph{(ii)} Write
\[
    Q_h(g)=\frac{1}{h}\int_0^{1-h}\|g(\xi+h)-g(\xi)\|_H\,d\xi.
\]
If $V_H(g;I)<\infty$, the variation-measure increment bound (cf.\ \cite[Remark~2.2, eq.~(2.3)]{Ambrosio1990metricBV}) and Tonelli's theorem give
\[
    Q_h(g)\le\frac{1}{h}\int_0^{1-h}|Dg|((\xi,\xi+h])\,d\xi
    \le |Dg|(I)=\Vxi(w).
\]
This proves \eqref{eqn:V-xi-translation}, the case $\Vxi(w)=\infty$ being immediate.

Conversely, suppose $L:=\sup_{0<h<1}Q_h(g)<\infty$. The Steklov averages
\[
    g_h(\xi)=\frac{1}{h}\int_\xi^{\xi+h}g(s)\,ds
\]
converge to $g$ in $L^1_{\mathrm{loc}}((0,1);H)$, and $g_h'(\xi)=[g(\xi+h)-g(\xi)]/h$ almost everywhere on $(0,1-h)$. Fix an open interval $K$ whose closure is a compact subset of $(0,1)$; for sufficiently small $h$, the averages $g_h$ are defined on $K$. Let $P_m$ be the orthogonal projection onto the span of the first $m$ vectors of an orthonormal basis of $H$. Since $P_mg_h\to P_mg$ in $L^1(K;H)$, finite-dimensional lower semicontinuity \cite[Remark~3.5 and Proposition~3.6]{AmbrosioFuscoPallara2000} and the contractivity of $P_m$ yield
\[
    V_H(P_mg;K)\le\liminf_{h\downarrow0}\int_K\|P_mg_h'(\xi)\|_H\,d\xi
    \le\liminf_{h\downarrow0}Q_h(g)\le L.
\]
To pass to $H$, note that $P_mv\to v$ in $H$ for every $v\in H$. Thus each finite partition sum for $g$, taken at its Lebesgue points in $K$, is the limit of the corresponding sums for $P_mg$. Taking the supremum over such partitions, and using contractivity for the opposite inequality, gives
\[
    V_H(g;K)=\sup_m V_H(P_mg;K)\le L.
\]
Finally, choosing $K=(\varepsilon,1-\varepsilon)$ and letting $\varepsilon\downarrow0$ yields $\Vxi(w)=\sup_{0<\varepsilon<1/2}V_H(g;(\varepsilon,1-\varepsilon))\le L$. Together with $Q_h(g)\le\Vxi(w)$, this proves \eqref{eqn:V-xi-sup-char}.
\end{proof}

\section{Existence of Optimal Controls: The Finite Agent Problem}\label{sec:existence}

This section establishes the existence of an optimal control pair for the $N$-agent problem, providing the starting point for the analysis of the full graphon control problem. We first show that each admissible control determines a unique state trajectory, and then prove that the cost $J_N$ attains its minimum over these control--state pairs.

Define for $\vec{x}_N \in \Rbb^N$, and $w \in \UbvN$
\begin{equation}\label{eqn: FN definition}
    F_N(\vec{x}_N, w) : = \frac{1}{N}\sum_{j=1}^N w_{i,j} \phi(x_j - x_i).
\end{equation} 
\begin{theorem}\label{thm:N-agent}
    Assume that $\phi \in \textrm{Lip}(\Rbb)$. Then for any $w:[0,T] \rightarrow \UbvN$ measurable, there exists a solution $\vec{x}_N \in AC([0,T] ; \Rbb^N)$ to the equation 
    \begin{equation}\label{eqn: N dynamics}
        \partial_t \vec{x}_N(t) = F_N(\vec{x}_N(t), w(t)), \qquad \vec{x}_N(0) = \vec{x}_{N,0} \in \Rbb^N
    \end{equation}
    where $F_N$ is defined by \eqref{eqn: FN definition}. Moreover, there exists a minimizing pair $\vec{x}_N^\star \in AC([0,T]; \Rbb^N)$ and $w^\star\in L^\infty([0,T]; \UbvN)$ of the functional $J_N$ defined in \eqref{eqn: JN definition} such that $(\vec{x}_N^\star, w^\star)$ that solves \eqref{eqn: N dynamics}.
\end{theorem}

\begin{proof} We prove existence of solutions for fixed graphons and existence of optimal controls separately. 

\medskip

   \noindent \textbf{Existence of Solutions:}
    Fix a measurable control $w : [0,T] \to \UbvN$. We establish existence and uniqueness of an absolutely continuous solution to \eqref{eqn: N dynamics} via a Picard fixed point argument in a time-weighted norm.
    
    Since $w(t) \in \UbvN \subset [0,\Lambda]^{N\times N}$ for a.e.\ $t$ and $\phi$ is Lipschitz, for any continuous curve $\vec{x}_N: [0,T] \to \Rbb^N$ the map $s \mapsto F_N(\vec{x}_N(s), w(s))$ is measurable and bounded. By Lebesgue's differentiation theorem, $\vec{x}_N \in AC([0,T]; \Rbb^N)$ solves \eqref{eqn: N dynamics} with $\vec{x}_N(0) = \vec{x}_{N,0}$ if and only if
    \begin{equation}\label{eqn: integral form}
        \vec{x}_N(t) = \vec{x}_{N,0} + \int_0^t F_N(\vec{x}_N(s), w(s))\,ds, \qquad t \in [0,T].
    \end{equation}
    
    Let $L_\phi$ denote the Lipschitz constant of $\phi$. We claim that
    \begin{equation}\label{eqn: FN Lipschitz}
        |F_N(\vec{x}_N, w) - F_N(\vec{y}_N, w)| \le L\, |\vec{x}_N - \vec{y}_N| \qquad \text{for all } \vec{x}_N, \vec{y}_N \in \Rbb^N,\ w \in \UbvN,
    \end{equation}
    with $L := L_\phi \sqrt{2(1+\Lambda)}$, independent of both $w$ and $N$. Set $p_{ij} := w_{ij}/N$; by the normalization \eqref{eqn: U nonsymmetric}, each row $(p_{ij})_{j=1}^N$ is a probability vector with $0 \le p_{ij} \le \Lambda/N$, and the column sums obey $\sum_{i} p_{ij} = \frac{1}{N} \sum_i w_{ij} \le \frac{1}{N} \cdot N \Lambda = \Lambda$. Setting $a_{ij} := \phi(x_j - x_i) - \phi(y_j - y_i)$, the elementary inequality $$a_{ij}^2 \le 2 L_\phi^2 (|x_i - y_i|^2 + |x_j - y_j|^2)$$ together with Jensen's inequality applied to the probability vector $(p_{ij})_{j=1}^N$ gives
    \[
        |F_N(\vec{x}_N,w) - F_N(\vec{y}_N,w)|^2
        = \sum_{i=1}^N \Big| \sum_{j=1}^N p_{ij} a_{ij} \Big|^2
        \le \sum_{i,j=1}^N p_{ij}\, a_{ij}^2.
    \]
    Using the row sums $\sum_j p_{ij} = 1$ and the column bound $\sum_i p_{ij} \le \Lambda$,
    \[
    \begin{aligned}
        2 L_\phi^2 \sum_{i,j} p_{ij}\bigl(|x_i - y_i|^2 + |x_j - y_j|^2\bigr) &= 2 L_\phi^2 \Big( \sum_i |x_i - y_i|^2 \sum_j p_{ij} + \sum_j |x_j - y_j|^2 \sum_i p_{ij} \Big) \\
        & \le 2 L_\phi^2 (1 + \Lambda)\,|\vec{x}_N - \vec{y}_N|^2,
    \end{aligned}
    \]
    which proves \eqref{eqn: FN Lipschitz}.
    
    Define the Picard operator
    \[
        \Phi : \Ccal([0,T]; \Rbb^N) \to \Ccal([0,T]; \Rbb^N),
        \qquad
        \Phi(\vec{x}_N)(t) := \vec{x}_{N,0} + \int_0^t F_N(\vec{x}_N(s), w(s))\,ds,
    \]
    which is well-defined by the measurability and boundedness noted above; its fixed points are exactly the solutions of \eqref{eqn: integral form}. For $\lambda > 0$ to be chosen, equip $\Ccal([0,T]; \Rbb^N)$ with the Bielecki norm
    \[
        \|\vec{x}_N\|_\lambda := \sup_{t \in [0,T]} e^{-\lambda t}\,|\vec{x}_N(t)|.
    \]
    This norm is equivalent to the uniform norm, so $(\Ccal([0,T]; \Rbb^N), \|\cdot\|_\lambda)$ is a Banach space. For $\vec{x}_N, \vec{y}_N \in \Ccal([0,T];\Rbb^N)$ and $t \in [0,T]$, \eqref{eqn: FN Lipschitz} gives
    \[
        e^{-\lambda t}\,|\Phi(\vec{x}_N)(t) - \Phi(\vec{y}_N)(t)|
        \le L \int_0^t e^{-\lambda(t-s)} \bigl(e^{-\lambda s}\, |\vec{x}_N(s) - \vec{y}_N(s)|\bigr)\, ds
        \le \frac{L}{\lambda}\,\|\vec{x}_N - \vec{y}_N\|_\lambda.
    \]
    Taking the supremum over $t$ yields
    \[
        \|\Phi(\vec{x}_N) - \Phi(\vec{y}_N)\|_\lambda \le \frac{L}{\lambda}\,\|\vec{x}_N - \vec{y}_N\|_\lambda.
    \]
    Choosing $\lambda>L$, the map $\Phi$ is a strict contraction on $(\Ccal([0,T];\Rbb^N), \|\cdot\|_\lambda)$.
    
    The Banach fixed point theorem produces a unique $\vec{x}_N^w \in \Ccal([0,T]; \Rbb^N)$ satisfying \eqref{eqn: integral form}. Since the right-hand side of \eqref{eqn: integral form} is the integral of an essentially bounded measurable function, $\vec{x}_N^w \in AC([0,T]; \Rbb^N)$ and satisfies \eqref{eqn: N dynamics} for a.e.\ $t \in [0,T]$. This completes the proof of existence and uniqueness of solutions for each fixed admissible control.

    \medskip
    
    \noindent \textbf{Existence of Optimal Control:} We apply Filippov's existence theorem for Bolza problems \cite[Theorem~9.3.i, p.~314]{Cesari1983}. The set $\UbvN$ is nonempty, compact, and convex: it contains $w\equiv1$, and its defining constraints consist of a bounded box, affine row-sum conditions, and a sublevel set of the continuous convex variation functional in \eqref{eqn:receiver-tv-discrete}. The preceding Lipschitz estimate and row normalization give
    \[
        |F_N(\vec{x}_N,w)| \le L|\vec{x}_N| + \sqrt{N}\,|\phi(0)|.
    \]
    By Gr\"onwall's inequality, every admissible trajectory is bounded uniformly
    \[
        \sup_{0\le t\le T}|\vec{x}_N(t)|
        \le \bigl(|\vec{x}_{N,0}| + T\sqrt{N}\,|\phi(0)|\bigr)e^{LT}.
    \]
    So we may restrict the time-state domain to a compact cylinder without excluding any admissible trajectory. Its product with $\UbvN$ is compact, the fixed initial condition defines a closed endpoint constraint, and the dynamics and running and terminal costs are continuous. The control $w\equiv1$, together with its solution constructed above, supplies an admissible pair.

    It remains to verify the theorem's convexity condition. Define the running cost
    \[
        \ell_N(\vec{x}_N,w)
        := \Var(\vec{x}_N) + \frac{\beta}{N^2}\sum_{i,j=1}^{N}|w_{i,j}|^2
    \]
    and, for fixed $\vec{x}_N$, the enlarged velocity-cost set
    \[
        \widetilde{Q}(\vec{x}_N)
        := \left\{(r,z)\in\Rbb^{N+1}:
        r\ge\ell_N(\vec{x}_N,w),\quad z=F_N(\vec{x}_N,w),\quad w\in\UbvN\right\}.
    \]
    Take $(r_k,z_k)\in\widetilde{Q}(\vec{x}_N)$, with corresponding controls $w_k$, for $k=1,2$. For $\theta\in[0,1]$, convexity of the admissible set gives $w_\theta=\theta w_1+(1-\theta)w_2\in\UbvN$. The dynamics are linear in $w$, while the discrete quadratic regularizer is convex for $\beta\ge0$. Since the variance term is independent of $w$, we obtain
    $$
    F_N(\vec{x}_N,w_\theta) = \theta z_1+(1-\theta)z_2,
    $$
    and
    $$
    \ell_N(\vec{x}_N,w_\theta) \le \theta\ell_N(\vec{x}_N,w_1)+(1-\theta)\ell_N(\vec{x}_N,w_2) \le \theta r_1+(1-\theta)r_2.
    $$
    Hence $\theta(r_1,z_1)+(1-\theta)(r_2,z_2)\in\widetilde{Q}(\vec{x}_N)$, proving convexity. All hypotheses of the existence theorem are satisfied, so $J_N$ attains its minimum at an admissible pair $(\vec{x}_N^\star,w^\star)$. Since the admissible controls are bounded, $w^\star\in L^\infty([0,T];\UbvN)$.
\end{proof}

\begin{remark}\label{rem:non-uniqueness}
Without the regularizing term, i.e. $\beta = 0$ in \eqref{eqn: JN definition}, minimizing controls are not unique in general. The following $N = 4$ example illustrates this. Let $\phi(z) = z$, $\xtar = 0$, and initial states $\vec{x}_4(0) = (1, -1, 1, -1)$, consider the two constant-in-time controls
    \begin{equation*}
        w_A = \begin{bmatrix}
            0 & 4 & 0 & 0 \\
            4 & 0 & 0 & 0 \\
            0 & 0 & 0 & 4 \\
            0 & 0 & 4 & 0
        \end{bmatrix},
        \qquad
        w_B = \begin{bmatrix}
            0 & 0 & 0 & 4 \\
            0 & 0 & 4 & 0 \\
            0 & 4 & 0 & 0 \\
            4 & 0 & 0 & 0
        \end{bmatrix}.
    \end{equation*}
Under $w_A$, agent $1, 2$ and $3, 4$ are paired; under $w_B$, agent $1, 4$ and $2, 3$ are paired. In either case, every agent is paired with an agent holding the opposite state, and solving  \eqref{eqn: N dynamics} yields the same trajectory $\vec{x}_4(t) = e^{-2t}(1, -1, 1, -1)$. 
\end{remark}

\section{Graphon Control for Infinite Agents}

\begin{theorem}\label{thm:graphon-existence}
    Let $\phi \in \text{Lip}(\Rbb)$, $x_{0} \in L^2(I)$, and $w$ be an admissible control. Then there exists unique $\mathcal{C}^1$ solutions $x: [0,T] \rightarrow L^2(I)$ to \eqref{eqn:graphon-dynamics} for all $T >1$. 
\end{theorem}

\begin{proof}
    The proof of \Cref{thm:N-agent} works with minor modifications on the Banach space $L^\infty( [0,T] ;L^2(I))$ with the norm
    \begin{equation}
        \Vert x \Vert_\lambda = \sup_{t\in[0,T]} e^{-\lambda t} \Vert x(t) \Vert_{L^2(I)},
    \end{equation}
    assuming we can prove that the map $F : L^2(I) \rightarrow L^2(I)$ defined by 
    \begin{equation}\label{eqn:F-graphon-def}
        F(x, w)(\xi) = \int_I w(\xi,\zeta) \phi( x(t,\zeta) - x(t,\xi)) d\zeta 
    \end{equation}
    is Lipchitz on $L^2(I)$. To do this, we simply check that for $x,y \in L^2(I)$ we have 
    \begin{align*}
        \Vert F(x,w) - F(y,w) \Vert^2 & = \int_I \left| \int_I w (\xi,\zeta) \left[ \phi(x(\zeta) - x(\xi) ) - \phi (y(\zeta) - y(\xi))\right] d\zeta \right|^2 d\xi \\
        & \leq \Vert w \Vert_{L^\infty}^2 [\phi]_{Lip}^2 \int_I \int_I \left| x(\zeta) - x(\xi)  - y(\zeta) + y(\xi) \right|^2 d\xi d\zeta \\
        & \leq 4 \Vert w \Vert_{L^\infty}^2 [\phi]_{Lip}^2\Vert x - y \Vert_{L^2{(I)}}^2
    \end{align*}
    Picking $\lambda = 8 \Vert w \Vert_{L^\infty}^2 [\phi]_{Lip}^2$ finishes the proof. Since every admissible control satisfies $\Vert w \Vert_{L^\infty(I^2)} \le \Lambda$ by \eqref{eqn:asymmetric-graphon}, one may take $\lambda = 8 \Lambda^2 [\phi]_{Lip}^2$, uniformly over the admissible class.
\end{proof}

We next establish that the projected $N$-agent problem is well posed, and satisfies some monotonicity properties. 

\begin{lemma} \label{lem:AN-monotone}
    Let $x(0) \in L^2(I)$ and fix a dyadic integer $M$. For any $N \geq M$ dyadic, the following minimization problem is well posed
    \begin{equation}\label{eqn:ANM}
    A_N(M) := \min \Big\{\mathcal{J}(w_N; x_M(0)): w_N \in L^\infty([0, T]; \UbvN)\Big\}
    \end{equation} 
    where $x_{M}(0) = P_M x(0)$ and $\mathcal{J}(\cdot, \cdot)$ is defined in \eqref{eqn:J-continuum-alt}. Moreover, we have that $A_N(M) \ge A_{2N}(M)$. 
\end{lemma}

\begin{figure}[ht]
    \centering
    \begin{subfigure}[c]{0.2\linewidth}
        \centering
        \includegraphics[width=\linewidth]{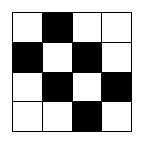}
        \caption{}
        \label{fig:graphon-N}
    \end{subfigure}%
    \begin{minipage}[c]{0.12\linewidth}
        \centering
        $\xrightarrow{\ \text{  Refine   }\ }$
    \end{minipage}%
    \begin{subfigure}[c]{0.2\linewidth}
        \centering
        \includegraphics[width=\linewidth]{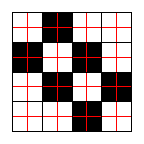}
        \caption{}
        \label{fig:graphon-2N}
    \end{subfigure}
    \caption{Refining the grid from $N \to 2N$.}
    \label{fig:graphons}
\end{figure}

\begin{proof}

Since $M\le N$ are dyadic, $P_Nx_M(0)=x_M(0)$. For histogram states
and kernels, the graphon equation on $I_i^N$ becomes
\[
    \partial_t  x_N^i(t)
    = \frac{1}{N}\sum_{j=1}^N
      \overline{w}_N^{i,j}(t)
      \phi\bigl(x_N^j(t)-x_N^i(t)\bigr),
\]
which is precisely the $N$-agent dynamics \eqref{eqn: N dynamics}. Moreover,
\[
    \overline{x}_N=\frac{1}{N}\sum_{i=1}^N x_N^i,
    \qquad
    \Var(x_N)=\frac{1}{N}\sum_{i=1}^N
      |x_N^i-\overline{x}_N|^2,
    \qquad
    \|w_N\|_{L^2(I^2)}^2
    =\frac{1}{N^2}\sum_{i,j=1}^N|\overline{w}_N^{i,j}|^2.
\]
Thus the histogram embedding preserves the dynamics and every term
of the cost. Conversely, uniqueness of \eqref{eqn:graphon-dynamics} as given by \Cref{thm:graphon-existence} identifies any solution with a
level-$N$ histogram control and initial datum $x_M(0)$ with this
embedding. \Cref{thm:N-agent} therefore supplies a minimizing pair
$(x_N^\star,w_N^\star)$, proving that $A_N$ is attained.

Now view these same functions on the refined $2N$-partition (refer to \Cref{fig:graphons}). Given the
    minimizer $(x_N^\star, w_N^\star)$, construct
    $(x_{2N}, w_{2N})$ at level $2N$ by splitting each cell
    $I_i^N$ into two halves $I_{2i-1}^{2N} \cup I_{2i}^{2N}$ and assigning the
    same nodal value (for the state) and the same block value (for the kernel)
    to both halves. By splitting $I_i^N = I_{2i-1}^{2N} \cup I_{2i}^{2N}$ and $I_i^N \times I_j^N = \cup_{a, b \in \{0, 1\}} I_{2i - a}^{2N} \times I_{2j-b}^{2N}$, we can immediately verify that 
    $$x_{2N}(t, \xi) = \sum_{\ell=1}^{2N}  x_{2N}^\ell \mathbf{1}_{I^{2N}_\ell} (\xi) = \sum_{i=1}^N x_N^{i\star}(t) \Big(\mathbf{1}_{I_{2i-1}^{2N}} + \mathbf{1}_{I_{2i}^{2N}}\Big) = \sum_{i=1}^N x_N^{i\star}(t) \mathbf{1}_{I_{i}^{N}}(\xi) = x_N^{\star}(t, \xi),
    $$
    $$
    w_{2N}(t, \xi, \zeta) = \sum_{p, q = 1}^{2N} \overline{w}_{2N}^{p, q} (t) \mathbf{1}_{I_{p}^{2N}}(\xi) \mathbf{1}_{I_{q}^{2N}}(\zeta) = \sum_{i, j = 1}^{N} \overline{w}_{N}^{i,j,\star}(t)  \mathbf{1}_{I_{i}^{N}}(\xi) \mathbf{1}_{I_{j}^{N}}(\zeta) = w_N^\star(t, \xi, \zeta).
    $$

Their histogram representations satisfy the box constraint, and
duplicating each column gives row sums $2N$.
\Cref{lem:variation-properties}(i) shows that their discrete variation
is unchanged, since both kernel representations have the same
continuum variation. Their initial datum remains
$x_M(0)=P_{2N}x_M(0)$. The pair is therefore admissible at level $2N$,
with unchanged dynamics and cost, so
\[
    A_{2N}\le J[x_N^\star,w_N^\star]=A_N.
\]

\end{proof}

We now prove some properties of the $P_N$ operators as defined in \Cref{def:PN}. In particular, they are $L^2$ projections, and are non-expansive in the bounded variation seminorm as defined in \Cref{def:U-BV}. 
\begin{lemma}\label{lem:PN-projection}
    The operators $P_N$ defined in \Cref{def:PN} are bounded operators from $L^2(I)$ to $L^2(I)$ independent of $N$. In particular, we have for any $x \in L^2(I)$ and any $w \in \Ubv$,
    \begin{equation}\label{eqn:PN-bounded}
        \Vert P_N x \Vert_{L^2(I)} \leq \Vert x \Vert_{L^2(I)}
    \end{equation}
    \begin{equation}\label{eqn:BV-bounded}
        \Vxi(P_N w) \le \Vxi(w)
    \end{equation}
\end{lemma}

\begin{proof}
By Cauchy--Schwarz on each cell,
\[
    \|P_Nx\|_{L^2(I)}^2
    = N\sum_{i=1}^N
      \left|\int_{I_i^N}x(\xi)\,d\xi\right|^2
    \le \|x\|_{L^2(I)}^2,
\]
which proves \eqref{eqn:PN-bounded}.

The variation bound is immediate for $N=1$. For $N\ge2$,
the cell averages in \eqref{eqn:PN-kernel} satisfy
\[
    \overline{w}_N^{i+1,j}-\overline{w}_N^{i,j}
    = N^2\int_{I_i^N}\int_{I_j^N}
      \bigl[w(\xi+1/N,\zeta)-w(\xi,\zeta)\bigr]
      \,d\zeta\,d\xi.
\]
Thus, applying \Cref{lem:variation-properties}(i), Cauchy--Schwarz
in the sender variable $\zeta$, and Minkowski's inequality give
\begin{align*}
    \Vxi(P_Nw)
    &= \sum_{i=1}^{N-1}
       \left(\frac{1}{N}\sum_{j=1}^N
       \left|\overline{w}_N^{i+1,j}
             -\overline{w}_N^{i,j}\right|^2\right)^{1/2} \\
    &\le N\int_0^{1-1/N}
       \|w(\xi+1/N,\cdot)-w(\xi,\cdot)\|_{L^2(I)}
       \,d\xi \\
    &\le \Vxi(w),
\end{align*}
where the last inequality follows from
\eqref{eqn:V-xi-translation} with $h=1/N$.
This proves \eqref{eqn:BV-bounded}.
\end{proof}

In the following lemma, we clarify in which topology the projections $P_N$ converge to the identity as $N\rightarrow \infty$. 

\begin{lemma}\label{lem:PN-strong}
Let $P_N$ be the projection operators of \Cref{def:PN}. Then for every $x \in L^2(I)$ and every $v \in L^2(I^2)$,
\begin{equation}\label{eqn:PN-strong}
    \Vert P_N x - x \Vert_{L^2(I)} \rightarrow 0
    \qquad \text{and} \qquad
    \Vert P_N v - v \Vert_{L^2(I^2)} \rightarrow 0
\end{equation}
as dyadic $N \rightarrow \infty$. Consequently, if $w \in L^2([0,T]; L^2(I^2))$, then $P_N w \rightarrow w$ in $L^2([0,T]; L^2(I^2))$, where $P_N$ acts on each time slice $w(t)$.
\end{lemma}
\begin{proof}
We give the argument on $I^2$; the argument on $I$ is identical. First let $g \in C(I^2)$. On each cell $I_i^N \times I_j^N$ the value of $P_N g$ is an average of values of $g$ over that cell, so
$$
\Vert P_N g - g \Vert_{L^\infty(I^2)} \le \omega_g\big(\sqrt{2}/N\big),
$$
where $\omega_g$ denotes the modulus of continuity of $g$. Since $g$ is uniformly continuous on the compact set $I^2$, we have $\omega_g(\sqrt{2}/N) \rightarrow 0$, so $P_N g \rightarrow g$ uniformly, and in particular in $L^2(I^2)$. Now let $v \in L^2(I^2)$ and fix $\varepsilon > 0$. By density, choose $g \in C(I^2)$ with $\Vert v - g \Vert_{L^2(I^2)} < \varepsilon$. Since \eqref{eqn:PN-bounded} holds verbatim on $L^2(I^2)$ (by the same cell-by-cell Cauchy--Schwarz computation, applied to cells $I_i^N \times I_j^N$ of area $1/N^2$),
$$
\Vert P_N v - v \Vert_{L^2(I^2)} \le \Vert P_N (v - g) \Vert_{L^2(I^2)} + \Vert P_N g - g \Vert_{L^2(I^2)} + \Vert g - v \Vert_{L^2(I^2)} \le 2\varepsilon + \Vert P_N g - g \Vert_{L^2(I^2)},
$$
whence $\limsup_{N} \Vert P_N v - v\Vert_{L^2(I^2)} \le 2\varepsilon$ for every $\varepsilon > 0$, proving \eqref{eqn:PN-strong}. Finally, for $w \in L^2([0,T]; L^2(I^2))$, the above gives $\Vert P_N w(t) - w(t) \Vert^2_{L^2(I^2)} \rightarrow 0$ for a.e.\ $t \in [0,T]$, with the domination $\Vert P_N w(t) - w(t) \Vert^2_{L^2(I^2)} \le 4 \Vert w(t) \Vert^2_{L^2(I^2)} \in L^1([0,T])$, so the time-integrated convergence follows from the dominated convergence theorem.
\end{proof}

In the following lemma, we establish that the $N$-agent optimal pairs $(x_N^\star, w_N^\star)$ as given in \Cref{lem:AN-monotone} have weak sub-sequential limits. We utilize apriori bounds on the solutions in addition to the monotonicity of \Cref{lem:AN-monotone} to establish this. 

\begin{lemma}\label{lem:weak-limits}
     Let $x(0) \in L^2(I)$ and fix dyadic $M \in \Nbb$. For each dyadic $N \geq M$, let $(x^\star_N, w^\star_N)$ be the optimal pair stated in \Cref{lem:AN-monotone} with initial data $x_M(0) = P_M x(0)$. Then there exist a subsequence (not relabelled) and limits $x^\star \in L^2([0, T]; L^2(I))$ and $w^\star \in L^2([0, T]; L^2(I^2))$ such that \begin{equation*}
         x^\star_N \rightharpoonup x^\star \in L^2([0,T]; L^2(I)), \quad \text{and} \quad w^\star_N \rightharpoonup w^\star \in L^2([0, T]; L^2(I^2)).
     \end{equation*}
     Moreover, $w^\star(t) \in \Ubv$ for a.e.\ $t$, i.e.\ $w^\star$ is an admissible control.
\end{lemma}

\begin{proof}
By \Cref{lem:AN-monotone},
$A_N=J[x_N^\star,w_N^\star]\le A_M<\infty$.
The pointwise constraint $0\le w_N^\star\le\Lambda$ also gives
$\|w_N^\star\|_{L^\infty([0,T]\times I^2)}\le\Lambda$,
so the controls are uniformly bounded in
$L^2([0,T];L^2(I^2))$.

To bound the states, since $\phi \in \text{Lip}(\Rbb)$, we can choose $C>0$ such that
$|\phi(z)|\le C(1+|z|)$. The dynamics, the bound
$w_N^\star\le\Lambda$, and Cauchy--Schwarz give, for a.e.\ $t$,
\begin{align*}
    \frac12\frac{d}{dt}\|x_N^\star(t)\|_{L^2(I)}^2
    &= \iint_{I^2} x_N^\star(t,\xi)w_N^\star(t,\xi,\zeta)
       \phi\bigl(x_N^\star(t,\zeta)-x_N^\star(t,\xi)\bigr)
       \,d\zeta\,d\xi\\
    &\le C\Lambda\Bigl(
        \|x_N^\star(t)\|_{L^1(I)}
        +\|x_N^\star(t)\|_{L^1(I)}^2
        +\|x_N^\star(t)\|_{L^2(I)}^2\Bigr)\\
    &\le C\Lambda\|x_N^\star(t)\|_{L^2(I)}
        +2C\Lambda\|x_N^\star(t)\|_{L^2(I)}^2.
\end{align*}

To handle a possibly vanishing norm, we record the regularization
argument used below. Suppose $\rho\ge0$, $\rho^2$ is absolutely
continuous, and
\begin{equation}\label{eqn:gronwall-regularized}
    \frac12\frac{d}{dt}\rho(t)^2
    \le a\rho(t)^2+b(t)\rho(t)
    \qquad\text{for a.e.\ }t\in[0,T],
\end{equation}
where $a>0$ and $b$ is nonnegative and locally integrable.
For $\varepsilon>0$, the function
$\rho_\varepsilon=(\rho^2+\varepsilon^2)^{1/2}$ is absolutely
continuous and satisfies
\[
    \rho_\varepsilon\rho_\varepsilon'
    \le a\rho_\varepsilon^2+b(t)\rho_\varepsilon.
\]
Dividing by $\rho_\varepsilon\ge\varepsilon$, applying
Gr\"onwall's inequality, and letting $\varepsilon\downarrow0$ yields
\begin{equation}\label{eqn:gronwall-conclusion}
    \rho(t)\le e^{at}\rho(0)
       +\int_0^t e^{a(t-s)}b(s)\,ds,
    \qquad t\in[0,T].
\end{equation}
Apply this with $\rho(t)=\|x_N^\star(t)\|_{L^2(I)}$,
$a=2C\Lambda$, and $b\equiv C\Lambda$.
Since the integral equals
$\tfrac12(e^{2C\Lambda t}-1)$, we obtain
\[
    \|x_N^\star(t)\|_{L^2(I)}
    \le e^{At}\bigl(\|x_M(0)\|_{L^2(I)}+B\bigr),
    \qquad A:=2C\Lambda,\quad B:=\tfrac12.
\]
Using $x_M(0)=P_Mx(0)$ and \eqref{eqn:PN-bounded}, this gives
\begin{equation}
    \|x_N^\star(t)\|_{L^2(I)}
    \le e^{At}\bigl(\|x(0)\|_{L^2(I)}+B\bigr).
\end{equation}
Thus the states are uniformly bounded in
$L^\infty([0,T];L^2(I))$. Reflexivity and the Banach--Alaoglu
theorem now yield a common subsequence, not relabelled, such that
\[
    x_N^\star\rightharpoonup x^\star
       \quad\text{in }L^2([0,T];L^2(I)),
    \qquad
    w_N^\star\rightharpoonup w^\star
       \quad\text{in }L^2([0,T];L^2(I^2)).
\]

We next verify admissibility of the limiting control. Define
\begin{equation}\label{eqn:K-BV-def}
\begin{aligned}
    K&:=\{w\in L^2([0,T]\times I^2):
              0\le w\le\Lambda\text{ a.e.}\},\\
    K_{BV}&:=
       \{w \in K :\Vxi(w(t))\le R\text{ for a.e. }t\in[0,T]\}.
\end{aligned}
\end{equation}
By \Cref{lem:variation-properties}(ii),
$\Vxi(v)=\sup_{0<h<1}F_h(v)$, where
\[
    F_h(v):=\frac1h\int_0^{1-h}
       \|v(\xi+h,\cdot)-v(\xi,\cdot)\|_{L^2(I_\zeta)}
       \,d\xi.
\]
Each $F_h$ is convex, so $K_{BV}$ is convex.
Moreover, the triangle inequality and Cauchy--Schwarz give
\[
    |F_h(u)-F_h(v)|
    \le\frac2h\int_I
       \|u(\xi,\cdot)-v(\xi,\cdot)\|_{L^2(I_\zeta)}\,d\xi
    \le\frac2h\|u-v\|_{L^2(I^2)}.
\]
If $w_n\to w$ strongly in $L^2([0,T]\times I^2)$ with
$w_n\in K_{BV}$, a subsequence satisfies
$w_n(t)\to w(t)$ in $L^2(I^2)$ for a.e.\ $t$.
The box constraint passes to the limit. On a common full-measure
set of times, the preceding estimate gives, for every $h\in(0,1)$,
\[
    F_h(w(t))=\lim_{n\to\infty}F_h(w_n(t))\le R.
\]
Taking the supremum over $h$ yields $\Vxi(w(t))\le R$.
Hence $K_{BV}$ is strongly closed and convex, and therefore
weakly closed by Mazur's theorem.
By \Cref{lem:variation-properties}(i), the discrete variation
constraint agrees with the continuum constraint on histogram
kernels, so $w_N^\star\in K_{BV}$ for every dyadic $N$ and consequently
$w^\star\in K_{BV}$ by weak convergence.

Finally, row normalization gives, for $\xi\in I_i^N$,
\[
    \int_I w_N^\star(t,\xi,\zeta)\,d\zeta
    =\frac1N\sum_{j=1}^N\overline{w}_N^{i,j}(t)=1.
\]
For measurable $E\subset I$ and $\eta\in L^2(0,T)$,
testing the weak convergence against
$\eta(t)\mathbf{1}_E(\xi)$ yields
\[
\begin{aligned}
    \int_0^T\eta(t)\int_E\int_I
       w^\star(t,\xi,\zeta)\,d\zeta\,d\xi\,dt
    &=\lim_{N\to\infty}
       \int_0^T\eta(t)\int_E\int_I
       w_N^\star(t,\xi,\zeta)\,d\zeta\,d\xi\,dt\\
    &=|E|\int_0^T\eta(t)\,dt.
\end{aligned}
\]
Since $E$ and $\eta$ are arbitrary,
$\int_Iw^\star(t,\xi,\zeta)\,d\zeta=1$ for a.e.\ $(t,\xi)$.
Together with the box and variation constraints, this proves
$w^\star(t)\in\Ubv$ for a.e.\ $t$, as required.
\end{proof}

The following lemma demonstrates that the cost functional is weakly lower semi-continuous along optimal pairs.
\begin{lemma}\label{lem:cost-lsc}
    Let $x(0) \in L^2(I)$ be fixed, and let $(x^\star_N, w^\star_N)$ be the optimal pairs with weak limits $(x^\star, w^\star)$ starting at initial data $x_M(0) = P_M x(0)$ for $N \geq M$ dyadic numbers. Then
    $$
    \liminf_{N \to \infty} J[x^\star_N, w^\star_N] \ge  J[x^\star, w^\star].
    $$
\end{lemma}

\begin{proof}
The integrated variance and squared control norm are convex
continuous functionals on their respective space-time $L^2$
spaces, and hence are weakly lower semicontinuous. Therefore,
\begin{equation} \label{eqn: running cost lsc}
\begin{aligned}
    \int_0^T\Var(x^\star(t))\,dt
    &\le\liminf_{N\to\infty}
       \int_0^T\Var(x_N^\star(t))\,dt,\\
    \beta\int_0^T\|w^\star(t)\|_{L^2(I^2)}^2\,dt
    &\le\liminf_{N\to\infty}
       \beta\int_0^T\|w_N^\star(t)\|_{L^2(I^2)}^2\,dt.
\end{aligned}
\end{equation}
The terminal cost requires additional time regularity.
Pass first to a subsequence realizing the liminf of the full cost,
without relabelling. By \Cref{lem:weak-limits},
\[
    \sup_{t\in[0,T]}\|x_N^\star(t)\|_{L^2(I)}
    \le e^{AT}\bigl(\|x(0)\|_{L^2(I)}+B\bigr)=:K.
\]
The integral form of the dynamics gives, for $0\le s\le t\le T$,
\begin{equation}\label{eqn:time-integral-deriv-norm}
    \|x_N^\star(t)-x_N^\star(s)\|_{L^2(I)}
    \le\int_s^t
       \|\partial_t x_N^\star(\tau)\|_{L^2(I)}\,d\tau.
\end{equation}
Using \eqref{eqn:graphon-dynamics},
$0\le w_N^\star\le\Lambda$, and
$|\phi(z)|\le C(1+|z|)$, Minkowski's inequality and
Cauchy--Schwarz give, for a.e.\ $\tau$,
\begin{equation}\label{eqn:unif-time-deriv-bound}
\begin{aligned}
    \|\partial_t x_N^\star(\tau)\|_{L^2(I)}
    \le C\Lambda
       \bigl(1+2\|x_N^\star(\tau)\|_{L^2(I)}\bigr)\le 2C\Lambda(1+K).
\end{aligned}
\end{equation}
Consequently, plugging in \eqref{eqn:unif-time-deriv-bound} into \eqref{eqn:time-integral-deriv-norm} gives
\begin{equation}\label{eqn:equicontinuity-in-time}
    \|x_N^\star(t)-x_N^\star(s)\|_{L^2(I)}
    \le 2C\Lambda(1+K)|t-s|,
\end{equation}
with constants independent of $N$.

Endow the ball
\[
    X:=\{x\in L^2(I):\|x\|_{L^2(I)}\le K\}
\]
with its weak topology, which is compact and metrizable since
$L^2(I)$ is separable and reflexive.
The preceding estimate gives equicontinuity in this topology,
so the metric Arzel\`a--Ascoli theorem yields a further subsequence
converging in $C([0,T];X)$.
Its limit agrees a.e.\ with the existing space-time weak limit
$x^\star$ and supplies its weakly continuous representative.
In particular,
\[
    x_N^\star(T)\rightharpoonup x^\star(T)
    \quad\text{in }L^2(I).
\]
Testing against $\mathbf{1}_I$ gives
$\overline{x}_N^\star(T)\to\overline{x}^\star(T)$, and hence
\[
    \gamma|\overline{x}_N^\star(T)-\xtar|^2
    \longrightarrow
    \gamma|\overline{x}^\star(T)-\xtar|^2.
\]
Combining this convergence with the two running-cost inequalities \eqref{eqn: running cost lsc}
along the chosen subsequence proves
\[
    J[x^\star,w^\star]
    \le\liminf_{N\to\infty}J[x_N^\star,w_N^\star].
\]
\end{proof}
\begin{remark}\label{rmk:unregularized-controls}
The dynamics \eqref{eqn:graphon-dynamics} provides no spatial smoothing, so under the pointwise constraints of $\mathcal{U}$ alone the optimal graphons may oscillate in $\xi$ at the grid scale, leaving no strongly convergent subsequence of states to pair against the weak convergence of the kernels. This variation bound on $w\in \Ubv$ and the bound $\phi \in L^\infty(\Rbb) \cap \textrm{Lip}(\Rbb)$ supplies the necessary regularity for $x\in L^2(I)$, as described in the following lemma.
\end{remark}

\begin{lemma}\label{lem:spatial-modulus}
Assume $\phi \in L^\infty(\Rbb) \cap \text{Lip}(\Rbb)$ and write $L_\phi := [\phi]_{\mathrm{Lip}}$. Let $x_0 \in L^2(I)$, let $w \in \mathcal{U}_{BV}$ be as in \Cref{def:U-BV}, and let $x$ be the solution of \eqref{eqn:graphon-dynamics} with initial datum $x_0$ given by \Cref{thm:graphon-existence}. For $h \in (0,1)$ define the spatial translation modulus
\begin{equation}\label{eqn:rho-h-def}
  \rho_h(t) := \bigl\| x(t,\cdot+h) - x(t,\cdot) \bigr\|_{L^2(I)},
\end{equation}
where the translation $\xi \mapsto \xi + h$ is taken modulo $1$. Then
\begin{equation}\label{eqn:spatial-modulus-bound}
  \sup_{t\in[0,T]} \rho_h(t)
  \;\le\;
  e^{L_\phi T}\Bigl( \bigl\| x_0(\cdot+h) - x_0 \bigr\|_{L^2(I)} + \frac{\|\phi\|_{L^\infty(\Rbb)}}{L_\phi}\,\sqrt{2 \Lambda (R + 2\Lambda)\,|h|} \Bigr).
\end{equation}
\end{lemma}

\begin{proof}
Fix $h \in (0,1)$ and set
\[
  \delta(t,\xi) := x(t,\xi+h) - x(t,\xi),
  \qquad
  \Delta_h w(\xi,\zeta) := w(\xi+h,\zeta) - w(\xi,\zeta),
\]
all translations taken on the torus, so that $\rho_h(t) = \|\delta(t,\cdot)\|_{L^2(I)}$. The constants below depend only on $\phi$, $\Lambda$, $R$, and $T$.

Evaluating \eqref{eqn:graphon-dynamics} at $\xi+h$ and at $\xi$ and subtracting,
add and subtract $w(\xi+h,\zeta)\,\phi(x(\zeta)-x(\xi))$ inside the integral to obtain
\begin{equation}\label{eqn:sm-splitting}
  \partial_t \delta(t,\xi) = A(t,\xi) + B(t,\xi),
\end{equation}
where
\begin{align*}
  A(t,\xi) &:= \int_I \Delta_h w(\xi,\zeta)\,\phi\bigl(x(t,\zeta)-x(t,\xi)\bigr)\,d\zeta,
  \\[2pt]
  B(t,\xi) &:= \int_I w(\xi+h,\zeta)\,
    \Bigl[\phi\bigl(x(t,\zeta)-x(t,\xi+h)\bigr) - \phi\bigl(x(t,\zeta)-x(t,\xi)\bigr)\Bigr]\,d\zeta.
\end{align*}
The term $A$ carries the difference of \emph{kernels}; the term $B$ carries the
difference of \emph{states}. Testing \eqref{eqn:sm-splitting} against $\delta$ and
integrating in $\xi$,
\begin{equation}\label{eqn:sm-energy}
  \tfrac12 \frac{d}{dt}\rho_h(t)^2
  = \int_I \delta(t,\xi) \,A(t,\xi) \, d\xi + \int_I \delta(t,\xi) \,B(t,\xi) \, d\xi.
\end{equation}

Since $\phi$ is $L_\phi$-Lipschitz and $\int_I w(\xi+h,\zeta)\,d\zeta = 1$, we obtain
\[
  |B(t,\xi)| \le L_\phi\,|\delta(t,\xi)| \int_I w(\xi+h,\zeta)\,d\zeta
             = L_\phi\,|\delta(t,\xi)|,
\]
so that
\begin{equation}\label{eqn:sm-B-bound}
  \int_I \delta(t,\xi) \,B(t,\xi) \,d\xi \le L_\phi \int_I |\delta(t,\xi)|^2\,d\xi = L_\phi\,\rho_h(t)^2.
\end{equation}
This term supplies the Gr\"onwall coefficient; it carries no smallness in $|h|$
and uses neither the $\mathrm{BV}$ bound nor boundedness of $\phi$.

Next, by Cauchy--Schwarz in the $\zeta$-integral,
\[
  |A(t,\xi)|^2
  \le \Bigl(\int_I |\Delta_h w(\xi,\zeta)|^2\,d\zeta\Bigr)
      \Bigl(\int_I \bigl|\phi(x(\zeta)-x(\xi))\bigr|^2\,d\zeta\Bigr).
\]
Writing $f_h(\xi) := \|\Delta_h w(\xi,\cdot)\|_{L^2(I_\zeta)} \ge 0$,
integrating in $\xi$ gives
\begin{equation}\label{eqn:sm-A-to-f}
  \|A(t,\cdot)\|_{L^2(I)}^2
  \le \|\phi\|_{L^\infty(\Rbb)}^2 \int_I f_h(\xi)^2\,d\xi.
\end{equation}
We will now interpolate $f_h$ in $L^\infty$ and $L^1$ to obtain quantitative smallness of $f$ in $L^2$. On one hand, the pointwise a.e.
bound $0 \le w \le \Lambda$ gives $\|w(\xi,\cdot)\|_{L^2(I_\zeta)} \le \Lambda$ for a.e.\ $\xi$, whence
\begin{equation}\label{eqn:sm-f-Linfty}
  \|f_h\|_{L^\infty(I)}
  \le \operatorname*{ess\,sup}_\xi\Bigl(\|w(\xi+h,\cdot)\|_{L^2(I_\zeta)} + \|w(\xi,\cdot)\|_{L^2(I_\zeta)}\Bigr)
  \le 2 \Lambda.
\end{equation}

Since the translation is taken modulo $1$ while the variation is measured on the interval, we split the integral: on $(0, 1-h)$ the torus translation coincides with the interval translation, so the vector-valued $\mathrm{BV}$ translation estimate (\Cref{lem:variation-properties}(ii), estimate \eqref{eqn:V-xi-translation}) applies, while on the strip $(1-h, 1)$, of measure $h$, we use the pointwise bound $f_h \le 2\Lambda$ of \eqref{eqn:sm-f-Linfty}. Together,
\begin{equation}\label{eqn:sm-f-L1}
  \|f_h\|_{L^1(I)}
  = \int_0^{1-h} f_h(\xi)\,d\xi + \int_{1-h}^{1} f_h(\xi)\,d\xi
  \le |h|\,\Vxi(w) + 2\Lambda\,|h|
  \le (R + 2\Lambda)\,|h|.
\end{equation}
Combining \eqref{eqn:sm-f-Linfty}--\eqref{eqn:sm-f-L1},
\[
  \int_I f_h(\xi)^2\,d\xi \le \|f_h\|_{L^\infty(I)}\,\|f_h\|_{L^1(I)} \le 2 \Lambda (R + 2\Lambda)\,|h|,
\]
so that \eqref{eqn:sm-A-to-f} becomes
\begin{equation}\label{eqn:sm-A-bound}
  \|A(t,\cdot)\|_{L^2(I)} \le \|\phi\|_{L^\infty(\Rbb)}\,\sqrt{2\Lambda (R + 2\Lambda)\,|h|}.
\end{equation}
By Cauchy--Schwarz in $\xi$,
$$
\int_I \delta\,A\,d\xi \le \rho_h(t)\,\|A(t,\cdot)\|_{L^2(I)}
 \le \rho_h(t)\,\|\phi\|_{L^\infty(\Rbb)}\sqrt{2\Lambda (R + 2\Lambda)\,|h|}.
 $$

Next, we use estimates to establish a Gr\"onwall inequality. Insert \eqref{eqn:sm-B-bound} and the bound just obtained into \eqref{eqn:sm-energy}
and abbreviate $S_h := \|\phi\|_{L^\infty(\Rbb)}\sqrt{2 \Lambda (R + 2\Lambda)\,|h|}$:
\begin{equation}\label{eqn:sm-diff-ineq}
  \tfrac12 \frac{d}{dt}\rho_h(t)^2 \le L_\phi\,\rho_h(t)^2 + S_h\,\rho_h(t).
\end{equation}
By the regularized Gr\"onwall argument recorded in the proof of \Cref{lem:weak-limits}, applied with $\rho = \rho_h$, $a = L_\phi$, and $b \equiv S_h$,
\[
  \rho_h(t)
  \le e^{L_\phi t}\rho_h(0)
     + \frac{e^{L_\phi t}-1}{L_\phi}\,S_h
  \le e^{L_\phi T}\Bigl( \rho_h(0) + \frac{S_h}{L_\phi} \Bigr),
  \qquad t \in [0,T],
\]
which is \eqref{eqn:spatial-modulus-bound} after taking the supremum over $t \in [0,T]$.
\end{proof}

\begin{corollary}
    Let $(x^\star_N, w^\star_N)_{N \ge M}$ be the optimal pairs of \Cref{lem:AN-monotone} with common initial datum $x_M(0) = P_M x(0)$, and let
    \begin{equation}
        \rho^N_h(t)  = \bigl\| x_N^\star(t,\cdot+h) - x_N^\star(t,\cdot) \bigr\|_{L^2(I)},
    \end{equation}
    with the translation taken modulo $1$.  Then
    \begin{equation}\label{eqn:spatial-modulus-uniform}
        \lim_{h \to 0} \sup_{N \geq M} \sup_{t \in [0, T]} \rho_h^N (t) = 0.
    \end{equation}
\end{corollary}

\begin{proof}
    Since $w^\star_N \in \mathcal{U}_{BV, N}$, so its histogram kernel lies in $\Ubv$by \Cref{lem:variation-properties}(i). The bound \eqref{eqn:spatial-modulus-bound} applies to every $x^\star_N$ with constants independent of $N$. Since $P_N P_M = P_M$ for dyadic $N \geq M$, so $x^\star_N(0) = P_N x_M(0) = x_M(0)$. The first term $\| x_M(0)(\cdot + h) - x_M(0) \|_{L^2(I)}$ of \eqref{eqn:spatial-modulus-bound} does not depend on $N$, and it vanishes as $h \to 0$ by the continuity of translation in $L^2(I)$, as we proved in \Cref{lem:PN-strong}. The second term of \eqref{eqn:spatial-modulus-bound} also vanishes as $h \to 0$, therefore gives \eqref{eqn:spatial-modulus-uniform}.
\end{proof}

\begin{remark} \label{rmk:spatial-modulus-general}
The boundedness of $\phi$ is not essential for \Cref{lem:spatial-modulus}; it enters only through \eqref{eqn:sm-A-to-f}. Alternatively, one could use the equation to show that $L^q$ moments of the state are propagated for $q > 2$ and run the same interpolation, obtaining the weaker modulus $\rho_h \lesssim |h|^{1/2 - 1/q}$, which still suffices for the compactness arguments below. We do not pursue this here.
\end{remark}

The next lemma is a quantitative stability estimate for the solution map of \eqref{eqn:graphon-dynamics} with respect to both the kernel and the initial datum. The lemma will also provide the passage $M \to \infty$ in the initial data in the proof of \Cref{thm:main}. Similar estimates can be found in \cite{doi:10.1137/130943741, Ayi2021}.

\begin{lemma}\label{lem:S-stability}
Assume $\phi \in L^\infty(\Rbb) \cap \textrm{Lip}(\Rbb)$ and write $L_\phi := [\phi]_{\mathrm{Lip}}$. Then for any two initial data $x_0, y_0 \in L^2(I)$ and any two admissible controls $v, w$, we have for all $t \in [0,T]$
\begin{equation}\label{eqn:S-stability}
    \Vert S(v; x_0)(t) - S(w; y_0)(t) \Vert_{L^2(I)} \le e^{2 L_\phi \Lambda t} \Big( \Vert x_0 - y_0 \Vert_{L^2(I)} + \Vert \phi \Vert_{L^\infty(\Rbb)} \int_0^t \Vert v(s) - w(s) \Vert_{L^2(I^2)} \, ds \Big).
\end{equation}
In particular, by the Cauchy--Schwarz inequality in time,
\begin{equation}\label{eqn:S-stability-C0}
    \sup_{t \in [0,T]} \Vert S(v; x_0)(t) - S(w; y_0)(t) \Vert_{L^2(I)} \le e^{2 L_\phi \Lambda T} \Big( \Vert x_0 - y_0 \Vert_{L^2(I)}
    + \sqrt{T} \, \Vert \phi \Vert_{L^\infty(\Rbb)} \, \Vert v - w \Vert_{L^2([0,T] \times I^2)} \Big).
\end{equation}
where $S(w; x_0)$ is the solution map of \eqref{eqn:graphon-dynamics}. 
\end{lemma}
\begin{proof}
Write $x := S(v; x_0)$, $y := S(w; y_0)$, and $\delta := x - y$, so that $\delta(0) = x_0 - y_0$. Subtracting the two equations and adding and subtracting $v(\xi,\zeta)\, \phi(y(\zeta) - y(\xi))$ inside the integral, and compute 
\[
\partial_t \delta(t,\xi) = A(t,\xi) + B(t,\xi) 
\]
where we are defining
\[
\begin{aligned}
A(t, \xi) &= \int_I \big( v(t,\xi,\zeta) - w(t,\xi,\zeta) \big)\, \phi\big(y(t,\zeta) - y(t,\xi)\big)\, d\zeta,  \\
B(t, \xi) &= \int_I v(t,\xi,\zeta) \Big[ \phi\big(x(t,\zeta) - x(t,\xi)\big) - \phi\big(y(t,\zeta) - y(t,\xi)\big) \Big] d\zeta.
\end{aligned}
\]
For the state-difference term $B$: since $\phi$ is $L_\phi$-Lipschitz,
\[
\begin{aligned}
|B(t,\xi)| &\le L_\phi \int_I v(t,\xi,\zeta)\, \big( |\delta(t,\zeta)| + |\delta(t,\xi)| \big)\, d\zeta\\
&\le L_\phi \Big( \Vert v(t,\xi,\cdot) \Vert_{L^2(I_\zeta)} \Vert \delta(t) \Vert_{L^2(I)} + |\delta(t,\xi)| \Big),
\end{aligned}
\]
where we used Cauchy--Schwarz on the first term and the row-stochasticity $\int_I v \, d\zeta = 1$ on the second. Since $0 \le v \le \Lambda$ and $\int_I v \, d\zeta = 1$, we have $\Vert v(t,\xi,\cdot)\Vert_{L^2(I_\zeta)}^2 \le \Lambda \int_I v \, d\zeta = \Lambda$, so taking $L^2$ norms in $\xi$,
\[
\begin{aligned}
\Vert B(t) \Vert_{L^2(I)} &\le L_\phi \big( \sqrt{\Lambda} + 1 \big) \Vert \delta(t) \Vert_{L^2(I)}\\
&\le 2 L_\phi \Lambda \Vert \delta(t) \Vert_{L^2(I)},
\end{aligned}
\]
using $1 + \sqrt{\Lambda} \le 2\Lambda$ for $\Lambda \ge 1$. For the kernel-difference term $A$: by Cauchy--Schwarz in $\zeta$ and the bound $|\phi| \le \Vert \phi \Vert_{L^\infty(\Rbb)}$,
$$
|A(t,\xi)| \le \Vert (v-w)(t,\xi,\cdot) \Vert_{L^2(I_\zeta)} \, \Vert \phi \Vert_{L^\infty(\Rbb)},
$$
hence
$$
\Vert A(t) \Vert_{L^2(I)} \le \Vert \phi \Vert_{L^\infty(\Rbb)} \Vert v(t) - w(t) \Vert_{L^2(I^2)}.
$$
Testing the equation for $\delta$ against $\delta$ and using Cauchy--Schwarz,
\[
\tfrac{1}{2} \tfrac{d}{dt} \Vert \delta(t) \Vert^2_{L^2(I)} \le 2 L_\phi \Lambda \Vert \delta(t) \Vert^2_{L^2(I)} + \Vert \phi \Vert_{L^\infty(\Rbb)} \Vert v(t) - w(t) \Vert_{L^2(I^2)} \Vert \delta(t) \Vert_{L^2(I)}.
\]
The regularized Gr\"onwall estimate \eqref{eqn:gronwall-conclusion}, applied with $\rho = \Vert \delta \Vert_{L^2(I)}$, $a = 2 L_\phi \Lambda$, and $b(t) = \Vert \phi \Vert_{L^\infty(\Rbb)} \Vert v(t) - w(t) \Vert_{L^2(I^2)}$ directly implies \eqref{eqn:S-stability}. The uniform bound \eqref{eqn:S-stability-C0} follows from 
$$
\int_0^t \Vert v - w \Vert_{L^2(I^2)}\, ds \le \sqrt{T} \Vert v - w \Vert_{L^2([0,T] \times I^2)}.
$$
\end{proof}

The next corollary shows the cost is Lipschitz in the initial data, which allows us pass from $x_M(0)$ to $x(0)$ in the proof of \Cref{thm:main}.

\begin{corollary}\label{cor:value-stability}
Assume $\phi \in L^\infty(\Rbb) \cap \textrm{Lip}(\Rbb)$. For every $\rho > 0$ there is a constant $C = C(\rho, T, \Lambda, \phi, \gamma)$ such that for all $x_0, y_0 \in L^2(I)$ with $\Vert x_0 \Vert_{L^2(I)}, \Vert y_0 \Vert_{L^2(I)} \le \rho$,
\begin{equation}\label{eqn:value-stability}
     \big| A_\infty(x_0) - A_\infty(y_0) \big| \leq \sup_{w \in L^\infty([0,T]; \Ubv)} \big| \mathcal{J}(w; x_0) - \mathcal{J}(w; y_0) \big| \le C \, \Vert x_0 - y_0 \Vert_{L^2(I)},
\end{equation}
where $\mathcal{J}(\cdot;\cdot)$ is defined in \eqref{eqn:J-continuum-alt} and $A_\infty(\cdot)$ is defined in \eqref{eqn:A-infty-y-def}.
\end{corollary}
\begin{proof}
{ Fix an admissible $w$ and write $x := S(w; x_0)$ and $y := S(w; y_0)$ for the trajectories issued from the two data. Arguing exactly as in the proof of \Cref{lem:weak-limits} we have the a priori bound
$$
\max\Big\{ \sup_{t \in [0,T]} \Vert x(t) \Vert_{L^2(I)}, \ \sup_{t \in [0,T]} \Vert y(t) \Vert_{L^2(I)} \Big\} \le e^{AT} (\rho + B) =: K_\rho,
$$
with $A, B$ depending only on $\phi$ and $\Lambda$. By the stability estimate \eqref{eqn:S-stability-C0} applied with the \emph{common} kernel $w$,
$$
\sup_{t \in [0,T]} \Vert x(t) - y(t) \Vert_{L^2(I)} \le e^{2 L_\phi \Lambda T} \, \Vert x_0 - y_0 \Vert_{L^2(I)}.
$$
The control penalty $\beta \Vert w \Vert^2_{L^2}$ is common to $\mathcal{J}(w; x_0)$ and $\mathcal{J}(w; y_0)$ and cancels in the difference, which therefore involves only the variance and terminal terms. Writing $Pu := u - \overline{u}$ for the (orthogonal, hence contractive) projection onto mean-zero functions, $\Var(u) = \Vert P u \Vert_{L^2(I)}^2$, so for $\Vert u \Vert, \Vert v \Vert \le K_\rho$
$$
|\Var(u) - \Var(v)| \le \big( \Vert Pu \Vert + \Vert Pv \Vert \big) \Vert P(u - v) \Vert \le 2 K_\rho \Vert u - v \Vert_{L^2(I)}.
$$
Since $|\overline{u} - \overline{v}| \le \Vert u - v \Vert_{L^2(I)}$ by Cauchy--Schwarz,
$$
\begin{aligned}
    \gamma \Big| |\overline{u} - \xtar|^2 - |\overline{v} - \xtar|^2 \Big| &= \gamma \, |\overline{u} - \overline{v}| \, |\overline{u} + \overline{v} - 2\xtar| \\
    &\le 2\gamma \big( K_\rho + |\xtar| \big) \Vert u - v \Vert_{L^2(I)}.
\end{aligned}
$$
Applying the last two displays to $u = x(t)$, $v = y(t)$ and combining with the stability estimate,
$$
\begin{aligned}
    \big| \mathcal{J}(w; x_0) - \mathcal{J}(w; y_0) \big| &\le \Big( 2 T K_\rho + 2 \gamma (K_\rho + |\xtar|) \Big) e^{2 L_\phi \Lambda T} \, \Vert x_0 - y_0 \Vert_{L^2(I)} \\
    &=: C \Vert x_0 - y_0 \Vert_{L^2(I)},
\end{aligned}
$$
uniformly over admissible $w$. The second inequality in \eqref{eqn:value-stability} follows from the elementary bound $|\inf_w f(w) - \inf_w g(w)| \le \sup_w |f(w) - g(w)|$.}
\end{proof}

The next lemma shows that the weak limits $(x^\star, w^\star)$ in \Cref{lem:weak-limits} are an optimal pair for the continuum problem with initial data $x_M(0)$. 

\begin{lemma}\label{lem:liminf}
Assume $\phi \in L^\infty(\Rbb) \cap \textrm{Lip}(\Rbb)$. Fix $x(0) \in L^2(I)$, a dyadic integer $M$, and set $x_M(0) := P_M x(0)$. Let $(x^\star_N, w^\star_N)_{N \ge M}$ be the optimal pairs constructed in \Cref{lem:AN-monotone} with initial datum $x_M(0)$, and let $(x^\star, w^\star)$ be a weak subsequential limits as constructed in \Cref{lem:weak-limits}. Then, the following conclusions hold.
\begin{enumerate}
    \item[(i)] Up to a subsequence, $x^\star_N \to x^\star$ strongly in $L^2([0,T] \times I)$;
    \item[(ii)] The pair $(x^\star, w^\star)$ is a solution to \eqref{eqn:graphon-dynamics} with $x^\star(0) = x_M(0)$;
    \item[(iii)] We have $\lim_{N \to \infty} A_N(M) = A_\infty(x_M(0))$, where the limit is over the full sequence. In particular the infimum in \eqref{eqn:A-infty-y-def} is attained, and $(x^\star, w^\star)$ is an optimal pair for the continuum problem with initial datum $x_M(0)$.
\end{enumerate}
\end{lemma}
\begin{proof}
Since $P_Nx_M(0)=x_M(0)$ for dyadic $N\ge M$, the
histogram--vector correspondence in the proof of
\Cref{lem:AN-monotone} gives
\begin{equation}\label{eqn:AN-kernel-form}
    A_N(M)
    =\min\bigl\{\mathcal{J}(w;x_M(0)):
        w\in L^\infty([0,T];\UbvN)\bigr\},
\end{equation}
where matrices are identified with their histogram kernels.
By \Cref{lem:variation-properties}(i), these kernels belong to
$\Ubv$, so $A_N(M)\ge A_\infty(x_M(0))$.

By \Cref{lem:weak-limits} and the compactness argument in the
proof of \Cref{lem:cost-lsc}, we may pass to a further subsequence,
without relabelling, such that
\[
    x_N^\star\to x^\star
    \quad\text{in }C([0,T];L^2_{\mathrm{weak}}(I)),
    \qquad
    \sup_N\sup_{t\in[0,T]}
       \|x_N^\star(t)\|_{L^2(I)}\le K.
\]
We use \eqref{eqn:spatial-modulus-uniform} to strengthen this
convergence. Let $\eta_\varepsilon$ be a standard mollifier on
the torus $I$, supported within distance $\varepsilon$ of zero.
Minkowski's inequality and the uniform spatial modulus give
\[
    \sup_{N\ge M}\sup_{t\in[0,T]}
    \|x_N^\star(t)*\eta_\varepsilon-x_N^\star(t)\|_{L^2(I)}
    \longrightarrow0
    \qquad\text{as }\varepsilon\downarrow0.
\]
For fixed $\varepsilon$, the mollified states are uniformly
bounded and uniformly Lipschitz in the spatial variable, by
Cauchy--Schwarz and the bound $K$. Arzel\`a--Ascoli therefore
makes the mollified family precompact in $L^2(I)$, and the
preceding approximation shows that $\{x_N^\star(t)\}_N$ is
precompact in $L^2(I)$ for every $t$. Since its weak limit is
$x^\star(t)$, we obtain
\[
    x_N^\star(t)\to x^\star(t)
    \quad\text{strongly in }L^2(I)
    \quad\text{for every }t\in[0,T].
\]
Weak lower semicontinuity gives
$\|x^\star(t)\|_{L^2(I)}\le K$. Thus
$\|x_N^\star(t)-x^\star(t)\|_{L^2(I)}^2\le4K^2$, and dominated
convergence yields
\[
    x_N^\star\to x^\star
    \quad\text{strongly in }L^2([0,T]\times I),
\]
which establishes point (i) of the lemma. 

We next identify the limiting dynamics. Write
$L_\phi:=[\phi]_{\mathrm{Lip}}$ and
$\Delta x(s,\xi,\zeta):=x(s,\zeta)-x(s,\xi)$.
For every $\eta\in L^2(I)$ and $t\in[0,T]$, the integral
formulation gives
\begin{equation}\label{eqn:weak-form-N}
\begin{aligned}
    \langle x_N^\star(t),\eta\rangle
    =\langle x_M(0),\eta\rangle+\int_0^t\iint_{I^2}
       w_N^\star(s,\xi,\zeta)
       \phi(\Delta x_N^\star(s,\xi,\zeta))\eta(\xi)
       \,d\zeta\,d\xi\,ds.
\end{aligned}
\end{equation}
The left-hand side converges to
$\langle x^\star(t),\eta\rangle$. For the integral, the bound
$w_N^\star\le\Lambda$, Lipschitz continuity of $\phi$, and
Cauchy--Schwarz give
\[
\begin{aligned}
    \left|\int_0^t\iint_{I^2}
       w_N^\star
       \bigl[\phi(\Delta x_N^\star)-\phi(\Delta x^\star)\bigr]
       \eta(\xi)\,d\zeta\,d\xi\,ds\right|\le
       2\Lambda L_\phi\sqrt{T}\,
       \|\eta\|_{L^2(I)}
       \|x_N^\star-x^\star\|_{L^2([0,T]\times I)}
       \longrightarrow0.
\end{aligned}
\]
Moreover, boundedness of $\phi$ ensures that
$\mathbf{1}_{[0,t]}(s)\phi(\Delta x^\star(s,\xi,\zeta))\eta(\xi)$
is a fixed element of $L^2([0,T]\times I^2)$, against which
$w_N^\star\rightharpoonup w^\star$ may be tested. Passing to
the limit in \eqref{eqn:weak-form-N}, we obtain
\[
    \langle x^\star(t),\eta\rangle
    =\langle x_M(0),\eta\rangle
      +\int_0^t
       \langle F(x^\star(s),w^\star(s)),\eta\rangle\,ds,
\]
with $F$ as in \eqref{eqn:F-graphon-def}. Row-stochasticity
gives
$\|F(x^\star(s),w^\star(s))\|_{L^2(I)}
 \le\|\phi\|_{L^\infty(\Rbb)}$,
so this identity is the integral equation in $L^2(I)$.
Uniqueness in \Cref{thm:graphon-existence} therefore yields that $x^\star, w^\star$ solves the equation \eqref{eqn:graphon-dynamics} with initial data $x_M(0)$. This establishes point (ii) of the lemma.

Finally, we work on proving point (iii) of the lemma. The control $w^\star$ is admissible by \Cref{lem:weak-limits}.
Hence
\begin{equation}\label{eqn:Ainfty-le-Jstar}
    A_\infty(x_M(0))\le J[x^\star,w^\star].
\end{equation}
By \Cref{lem:AN-monotone}, the nonnegative sequence $A_N(M)$
is non-increasing, so its limit exists over the full dyadic
sequence. Applying \Cref{lem:cost-lsc} along the chosen
subsequence gives
\begin{equation}\label{eqn:limAN-ge-Jstar}
    J[x^\star,w^\star]\le\lim_{N\to\infty}A_N(M).
\end{equation}

For the reverse inequality, fix
$w\in L^\infty([0,T];\Ubv)$ and set $w_N(t):=P_Nw(t)$. We will use the sequence $w_N$ as a recovery sequence for $A_\infty(x_M(0))$. Note, the bounded linear operator $P_N$ preserves measurability,
and cell averaging preserves the box and row-sum constraints.
Moreover, \Cref{lem:PN-projection} gives
\[
    \Vxi(w_N(t))\le\Vxi(w(t))\le R
    \qquad\text{for a.e. }t.
\]
Thus $w_N\in L^\infty([0,T];\UbvN)$, and
\eqref{eqn:AN-kernel-form} implies
\begin{equation}\label{eqn:AN-le-recovery}
    A_N(M)\le\mathcal{J}(w_N;x_M(0)).
\end{equation}
By \Cref{lem:PN-strong},
$w_N\to w$ strongly in $L^2([0,T]\times I^2)$.
The stability estimate \eqref{eqn:S-stability-C0} then gives
\[
    \sup_{t\in[0,T]}
    \|S(w_N;x_M(0))(t)-S(w;x_M(0))(t)\|_{L^2(I)}
    \longrightarrow0.
\]
The variance and terminal cost are continuous uniformly on
bounded subsets of $L^2(I)$, while strong convergence of the
controls gives convergence of their squared $L^2$ norms.
Consequently,
$\mathcal{J}(w_N;x_M(0))\to\mathcal{J}(w;x_M(0))$.
Passing to the limit in \eqref{eqn:AN-le-recovery} and taking
the infimum over admissible $w$ yields
\[
    \lim_{N\to\infty}A_N(M)\le A_\infty(x_M(0)).
\]
Together with \eqref{eqn:Ainfty-le-Jstar} and
\eqref{eqn:limAN-ge-Jstar}, this proves
\[
    A_\infty(x_M(0))
    =J[x^\star,w^\star]
    =\lim_{N\to\infty}A_N(M),
\]
so $(x^\star,w^\star)$ attains the infimum in
\eqref{eqn:A-infty-y-def}.
\end{proof}

Fix a finite time horizon $T>0$. The following theorem is the main result of this paper: it establishes the existence of continuum optimal controls, convergence of the discrete optimal values, and subsequential convergence of optimal pairs.

\begin{theorem}\label{thm:main}
Assume $\phi \in L^\infty(\Rbb) \cap \textrm{Lip}(\Rbb)$, and let $x(0) \in L^2(I)$. Consider the cost $J$ and admissible control classes introduced in \Cref{sec:setting}. Then the following hold.
\begin{enumerate}
    \renewcommand{\labelenumi}{(\roman{enumi})}
    \item \emph{(Existence.)} The continuum problem admits an optimal pair $(x^\star, w^\star)$, with $w^\star \in L^\infty([0,T]; \Ubv)$. The pair solves the graphon dynamics \eqref{eqn:graphon-dynamics} with initial condition $x^\star(0) = x(0)$ and minimizes $J$ among all such admissible pairs.
    \item \emph{(Convergence of optimal values.)} Let $A_N(M)$ be the minimum cost for the $N$-agent problem with initial datum $P_M x(0)$, as defined in \Cref{lem:AN-monotone}. Then
    $$
    \lim_{M \to \infty} \, \lim_{N \to \infty} A_N(M) = J[x^\star, w^\star].
    $$
    Both limits exist and are taken through dyadic integers, with $N \ge M$; for each fixed $M$, the sequence $A_N(M)$ is non-increasing.
    \item \emph{(Subsequential convergence of optimal pairs.)} Up to a subsequence of dyadic integers $M \to \infty$ (not relabelled), one can choose dyadic $N = N(M) \ge M$ and optimal pairs $(x^\star_N, w^\star_N)$ for the $N$-agent problems with initial datum $P_M x(0)$ such that, as $M \to \infty$,
    $$
    w^\star_N \rightharpoonup w^\star \ \text{ weakly in } L^2([0,T] \times I^2),
    \qquad
    x^\star_N \to x^\star \ \text{ strongly in } L^2([0,T] \times I),
    $$
    and $J[x^\star_N, w^\star_N] \to J[x^\star, w^\star]$, where the limiting pair $(x^\star, w^\star)$ is an optimal pair as in item~(i).
\end{enumerate}
\end{theorem}

\begin{remark}
{\color{black} No uniqueness of the optimal pair is asserted. In item~(iii), a suitable dyadic choice $N=N(M)\ge M$ yields convergence of optimal pairs up to a subsequence (not relabelled). The iterated limits of the optimal values in item~(ii) hold over the full dyadic sequences.}
\end{remark}

\begin{proof}
{\color{black}
We begin with the value convergence required for item~(ii). Fix a dyadic $M$ and set $x_M(0)=P_Mx(0)$. By \Cref{lem:AN-monotone}, each dyadic $N\ge M$ admits an optimal pair $(x_N^\star,w_N^\star)$ with this initial datum. By \Cref{lem:liminf},
$A_N(M)\downarrow A_\infty(x_M(0))$ as $N\to\infty$, and, up to a subsequence in $N$ (not relabelled),
\[
    w_N^\star\rightharpoonup w^\star(M)
        \quad\text{in }L^2([0,T]\times I^2),
    \qquad
    x_N^\star\to x^\star(M)
        \quad\text{in }L^2([0,T]\times I).
\]
The limiting pair is optimal for the continuum problem with initial datum $x_M(0)$, so
\begin{equation}\label{eqn:thm-main-step1}
    J[x^\star(M),w^\star(M)]=A_\infty(x_M(0)).
\end{equation}

By \Cref{lem:PN-strong}, $x_M(0)\to x(0)$ strongly in $L^2(I)$. The projection bound \eqref{eqn:PN-bounded} places these data in the ball of radius $\rho:=1+\|x(0)\|_{L^2(I)}>0$. Hence \Cref{cor:value-stability} gives
\begin{equation}\label{eqn:thm-main-step2}
    \big|A_\infty(x_M(0))-A_\infty(x(0))\big|
    \le C(\rho)\,\|x_M(0)-x(0)\|_{L^2(I)}\longrightarrow0.
\end{equation}
Together with the fixed-$M$ value convergence, this identifies the iterated limit in item~(ii) as $A_\infty(x(0))$. To prove item~(i) and complete item~(ii), we now show that this value is attained.

The controls $w^\star(M)$ belong to the set $K_{BV}$ in \eqref{eqn:K-BV-def}, which imposes the box and receiver-variation constraints. As established in the proof of \Cref{lem:weak-limits}, this set is bounded, convex and weakly closed. Thus, up to a subsequence in $M$ (not relabelled),
\[
    w^\star(M)\rightharpoonup w^\star
        \quad\text{in }L^2([0,T]\times I^2),
    \qquad w^\star\in K_{BV}.
\]
Row normalization also passes to the limit: the map $w\mapsto\int_Iw\,d\zeta$ is bounded and linear from $L^2([0,T]\times I^2)$ to $L^2([0,T]\times I)$, and hence is weakly continuous. Consequently $w^\star\in L^\infty([0,T];\Ubv)$, as in the admissibility argument at the end of the proof of \Cref{lem:weak-limits}.

To obtain compactness of the states, boundedness of $\phi$ and row normalization give
\[
    \|\partial_t x^\star(M)(t)\|_{L^2(I)}\le\|\phi\|_{L^\infty(\Rbb)}
    \qquad\text{for a.e. }t.
\]
Therefore
\begin{align*}
    \sup_M\sup_{t\in[0,T]}\|x^\star(M)(t)\|_{L^2(I)}
        &\le\|x(0)\|_{L^2(I)}+T\|\phi\|_{L^\infty(\Rbb)},\\
    \|x^\star(M)(t)-x^\star(M)(s)\|_{L^2(I)}
        &\le\|\phi\|_{L^\infty(\Rbb)}|t-s|.
\end{align*}
The state bound places these trajectories in a fixed weakly compact, metrizable ball of $L^2(I)$, and the time estimate gives equicontinuity in that topology. Arzel\`a--Ascoli therefore yields, up to a further subsequence (not relabelled),
\[
    x^\star(M)\to x^\star
        \quad\text{in }C([0,T];L^2_{\mathrm{weak}}(I)).
\]

Since $x_M(0)\to x(0)$ in $L^2(I)$, the family of initial data is precompact. Continuity of translations is therefore uniform on this family, giving
\[
    \sup_M\|x_M(0)(\cdot+h)-x_M(0)\|_{L^2(I)}\longrightarrow0
        \qquad\text{as }h\to0,
\]
where translations are taken modulo $1$. The spatial estimate of \Cref{lem:spatial-modulus}, using a positive Lipschitz bound for $\phi$, then gives the same uniform translation continuity for the states, uniformly in time. Together with the state bound, the Fr\'echet--Kolmogorov criterion gives precompactness in $L^2(I)$ at each time. The weak convergence already obtained therefore implies strong convergence at each time, and dominated convergence gives
\[
    x^\star(M)\to x^\star
        \quad\text{strongly in }L^2([0,T]\times I).
\]

For each $M$, the state $x^\star(M)$ satisfies
\[
    \partial_s x^\star(M)(s,\xi)
    =\int_I w^\star(M)(s,\xi,\zeta)\,
      \phi\bigl(x^\star(M)(s,\zeta)-x^\star(M)(s,\xi)\bigr)\,d\zeta
\]
for a.e.\ $(s,\xi)$, with initial datum $x_M(0)$. We pass to the limit in its integral formulation \eqref{eqn:weak-form-N}. Fix $t\in[0,T]$ and $\eta\in L^2(I)$, and compare the interaction integral for $(x^\star(M),w^\star(M))$ with that for $(x^\star,w^\star)$. Adding and subtracting the interaction formed with $w^\star(M)$ and $x^\star$ splits their difference into the sum
\[
\begin{aligned}
    &\int_0^t\iint_{I^2} w^\star(M)(s,\xi,\zeta)
      \Bigl[\phi\bigl(x^\star(M)(s,\zeta)-x^\star(M)(s,\xi)\bigr)
      -\phi\bigl(x^\star(s,\zeta)-x^\star(s,\xi)\bigr)\Bigr]
      \eta(\xi)\,d\zeta\,d\xi\,ds\\[0.75em]
    &{}+\int_0^t\iint_{I^2}\bigl[w^\star(M)(s,\xi,\zeta)-w^\star(s,\xi,\zeta)\bigr]
      \phi\bigl(x^\star(s,\zeta)-x^\star(s,\xi)\bigr)
      \eta(\xi)\,d\zeta\,d\xi\,ds.
\end{aligned}
\]
The first integral tends to zero by the state-difference estimate immediately following \eqref{eqn:weak-form-N} in the proof of \Cref{lem:liminf}(ii), and the second by the fixed-test-function argument in the same passage. Those arguments apply with $(x_N^\star,w_N^\star)$ replaced by $(x^\star(M),w^\star(M))$, using the strong state convergence and weak control convergence proved above. The initial-data term converges because $x_M(0)\to x(0)$ strongly in $L^2(I)$. Together with $x^\star(M)(t)\rightharpoonup x^\star(t)$, this proves that $(x^\star,w^\star)$ solves \eqref{eqn:graphon-dynamics} with $x^\star(0)=x(0)$.

Finally, the running costs are convex continuous quadratic functionals and hence weakly lower semicontinuous. Weak convergence of $x^\star(M)(T)$ to $x^\star(T)$ implies convergence of the terminal means by testing against $\mathbf{1}_I$. Combining these facts gives
\[
    J[x^\star,w^\star]
    \le\liminf_{M\to\infty}J[x^\star(M),w^\star(M)]
    =\lim_{M\to\infty}A_\infty(x_M(0))
    =A_\infty(x(0)),
\]
where the equalities follow from \eqref{eqn:thm-main-step1} and \eqref{eqn:thm-main-step2}. Since $(x^\star,w^\star)$ is admissible, the reverse inequality follows from the definition of $A_\infty(x(0))$. This proves item~(i). The identity $J[x^\star,w^\star]=A_\infty(x(0))$ also identifies the iterated limit established above, completing item~(ii).

It remains to choose the discrete optimal pairs in item~(iii). On the bounded set $K_{BV}$ defined in \eqref{eqn:K-BV-def}, the weak topology of the separable Hilbert space $L^2([0,T]\times I^2)$ is metrizable; fix a metric $d$ inducing this topology. Keep the subsequence of dyadic $M$ chosen above. For each such $M$, recall that $(x^\star(M),w^\star(M))$ is the continuum optimal pair in \eqref{eqn:thm-main-step1}, obtained as a subsequential limit of $N$-agent optimal pairs $(x_N^\star,w_N^\star)$ as $N\to\infty$ with $M$ fixed. Both pairs have the same initial datum $P_Mx(0)$. Taking $N$ sufficiently far along this fixed-$M$ subsequence, we can choose a dyadic $N=N(M)\ge M$ such that
\begin{align*}
    \big|A_N(M)-A_\infty(x_M(0))\big|&\le M^{-1},\\
    d\bigl(w_N^\star,w^\star(M)\bigr)&\le M^{-1},\\
    \|x_N^\star-x^\star(M)\|_{L^2([0,T]\times I)}&\le M^{-1}.
\end{align*}
Along this choice $N=N(M)$, the triangle inequality yields
\[
    w_N^\star\rightharpoonup w^\star
        \quad\text{in }L^2([0,T]\times I^2),
    \qquad
    x_N^\star\to x^\star
        \quad\text{in }L^2([0,T]\times I).
\]
The value estimate and \eqref{eqn:thm-main-step2} also give
\[
    J[x_N^\star,w_N^\star]=A_N(M)\longrightarrow
        A_\infty(x(0))=J[x^\star,w^\star].
\]
This proves item~(iii).}
\end{proof}

The next remark shows that increasing the receiver-variation bound recovers the optimal value of the problem without a variation constraint.
\begin{remark}\label{rmk:value-invariance}
Under the assumptions of \Cref{thm:main}, keep $\Lambda$ and $x(0)$ fixed and write $\Ubv(R):=\{w\in\mathcal{U}:\Vxi(w)\le R\}$ to make the radius in \eqref{eqn:U-BV-def} explicit. Then
\[
\inf_{w\in L^\infty([0,T];\Ubv(R))}\mathcal{J}(w;x(0))
\;\downarrow\;
\inf_{w\in L^\infty([0,T];\mathcal{U})}\mathcal{J}(w;x(0))
\qquad\text{as }R\to\infty.
\]
Indeed, fix any admissible $w$ taking values in $\mathcal{U}$ and extend it from $I=[0,1]$ to $\Rbb$ as a $1$-periodic function of the receiver variable $\xi$. For a nonnegative smooth periodic approximate identity $\rho_\epsilon$ of unit mass, set $w^\epsilon=\rho_\epsilon *_\xi w$. Mollification preserves measurability and the constraints of $\mathcal{U}$, and $w^\epsilon\to w$ strongly in $L^2([0,T]\times I^2)$. Since $\|w(t,\xi,\cdot)\|_{L^2(I)}\le\sqrt{\Lambda}$, we also have
\[
\Vxi(w^\epsilon(t))
\le\int_I\|\partial_\xi w^\epsilon(t,\xi,\cdot)\|_{L^2(I)}\,d\xi
\le\sqrt{\Lambda}\,\|\rho_\epsilon'\|_{L^1(I)}=:R_\epsilon
\]
for a.e.\ $t$, so $w^\epsilon\in L^\infty([0,T];\Ubv(R_\epsilon))$. The stability estimate \eqref{eqn:S-stability-C0}, whose proof uses no receiver-variation bound, gives convergence of the corresponding states in $C([0,T];L^2(I))$. All terms of the cost \eqref{eqn:J-continuum} therefore converge. The restricted infima decrease with $R$ and are bounded below by the unrestricted infimum; approximating each unrestricted competitor as above gives the reverse inequality for their limit.
\end{remark}

\nocite*
\bibliographystyle{plain}
\bibliography{references}

\newpage
\section{Appendix: Numerical Experiments}\label{sec:numerics}

In this appendix, we present two numerical experiments on a smooth bounded confidence opinion model, and a Kuramoto model. We design a neural network architecture and training protocol to learn an approximate optimal control, motivated by the $N$-agent constructions in \Cref{sec:existence}. In \Cref{sec:numerical-method}, we describe the discretizations scheme and network construction, then continue on to \Cref{subsec:exp-protocol} to describe our experimental protocol. \Cref{subsec:exp-bc} describes the phenomenon of clustering for a smooth bounded confidence opinion model with and without a control. \Cref{subsec:exp-kuramoto} describes untwisting twisted states of a Kuramoto model. Throughout these last two sections, we keep track of the learned receiver-variation bound $V_\xi(\widehat{w}(t))$.

\subsection{The method}\label{sec:numerical-method}

We consider the finite-agent dynamics, allowing additionally for stochasticity $\sigma \geq 0$,
$$
d x_N^i = \frac{1}{N} \sum_{j=1}^N \overline{w}_N^{i,j}(t) \, \phi(x_N^j - x_N^i) \, dt + \sigma \, dB_i^t, \qquad i \in \{1, \cdots, N\}, \quad t > 0,
$$
where $\{B_i^t\}_{i=1}^N$ are independent and identically distributed Brownian motions and $\overline{w}_N^{i,j}$ is the cell average \eqref{eqn: finite graphon approximation B}.  We shall be working in discrete time, with time step $\Delta t > 0$ and total number of time steps $N_T = T/\Delta t$. Applying explicit Euler results in the dynamics for each time step $n = 1, \dots, N_T$, with $s_n = n \Delta t$ and
\begin{equation}
    x_N^i(s_n) = x_N^i(s_{n-1}) + \frac{\Delta t}{N} \sum_{j=1}^N \overline{w}_N (s_{n-1}, \xi_i,\zeta_j) \phi(x_N^j(s_{n-1}) - x_N^i(s_{n-1})) + \sigma \Delta B_{n-1}^i.
\end{equation}
Here, $\Delta B_{n-1}^i \sim \mathcal{N}(0, \Delta t)$ are i.i.d samples from the Gaussian distribution with mean $0$ and variance $\Delta t$. Let $w = w_{\theta}$ be an approximation of the control via the deep neural network of \Cref{sec:nn-arch}. The parametrized cost is the cost $J_N$ of \eqref{eqn: JN definition} evaluated along the controlled trajectories. After discretizing time and simplifying (see \Cref{lem:AN-monotone}) we define
\begin{equation}\label{eqn:J-cost-1}
    \begin{aligned} 
        J^{N,\Delta t} (\theta) &= {\mathbb{E}}\Bigg[ \sum_{n=0}^{N_T-1}  \bigg(\frac{1}{N} \sum_{k=1}^N (x_N^k (s_n))^2 - \frac{1}{N^2} \sum_{i=1}^N \sum_{j=1}^N x_N^i (s_n) x_N^j (s_n) \bigg) \Delta t \\
        & \qquad + \beta \sum_{n=0}^{N_T-1}\|w_{\theta}(s_n)\|^2_{L^2(I^2)}  \Delta t 
        + \gamma | \overline{x}_N(T) - \xtar|^2 \Bigg],
    \end{aligned}
\end{equation}
where $\mathbb{E}$ is the expectation over the noise.
\subsubsection{Neural network architecture}\label{sec:nn-arch}

We now describe a neural parameterization of the interaction kernel in which the control depends on time and labels only (and \emph{not} on the current opinion state).
In particular, we parameterize a time-dependent graphon
\begin{equation}\label{eqn:nn-def-labelonly}
    w_\theta : [0,T]\times I^2 \to [0,\Lambda], \qquad (t,\xi,\zeta)\mapsto w_\theta(t,\xi,\zeta).
\end{equation}
The controlled dynamics remain state-dependent through the evolution operator (cf.\ \eqref{eqn:graphon-dynamics}), but the controller itself is \emph{label-only}. We impose the structural constraints
\begin{align}
    &0 \le w_\theta(t,\xi,\zeta)\le \Lambda, \label{eqn:nn-range-labelonly}\\
    &\int_I w_\theta(t,\xi,\zeta)\,d\zeta = 1 \quad \text{for all }(t,\xi)\in[0,T]\times I, \label{eqn:nn-rownorm-labelonly}
\end{align}
which encode a fixed attention budget per receiver label. We discretize $I=[0,1]$ using nodes $\{\xi_i\}_{i=1}^N\subset I$ and represent the kernel by an $N\times N$ matrix $W_\theta(t)$ with entries
\begin{equation}\label{eqn:nn-discrete-labelonly}
    W_{\theta,ij}(t) \approx \frac{1}{N} \overline{w}_N^{i,j}(t) = \frac{1}{N} \overline{w}_N(t, \xi_i, \xi_j)
\end{equation}
To represent $\xi\in I$ we use a Fourier feature map $\psi:I\to\Rbb^{2m}$, for example
\[
    \psi(\xi)=(\sin(2\pi 2^k\xi),\cos(2\pi 2^k\xi))_{k=0}^{m-1}.
\]
Define receiver and sender embeddings
\begin{equation}\label{eqn:nn-embed-labelonly}
    h_i := f_\theta\!\bigl(\psi(\xi_i)\bigr)\in\Rbb^d, \qquad \widetilde{h}_j := \widetilde{f}_\theta\!\bigl(\psi(\xi_j)\bigr)\in\Rbb^d
\end{equation}
where $f_\theta, \widetilde{f}_\theta: \Rbb^{2m} \rightarrow \Rbb^d$ are a feed-forward neural network.
We then define scores
\begin{equation}\label{eqn:nn-score-labelonly}
    s_{ij}(t) := g_\theta\!\bigl(t, h_i, \widetilde{h}_j\bigr)\in\Rbb,
\end{equation}
where $g_\theta:\Rbb^{1+ 2d}\to\Rbb$ is a feed-forward neural network. Finally, we obtain the weight matrix via a row-softmax:
\begin{equation}\label{eqn:nn-softmax-labelonly}
    W_{\theta,ij}(t) := \frac{\exp(s_{ij}(t))}{\sum_{k=1}^{N}\exp(s_{ik}(t))},
    \qquad 1\le i,j\le N.
\end{equation}
By construction, $W_{\theta,ij}(t) > 0$ and $\sum_{j=1}^N W_{\theta,ij}(t)=1$.

\subsection{Experimental protocol}\label{subsec:exp-protocol}

The $N$-agent dynamics \eqref{eqn: N dynamics} are integrated by the explicit Euler scheme of \Cref{sec:numerical-method}, and the control is parametrized by the row-softmax network of \eqref{eqn:nn-softmax-labelonly}. We trained the network by stochastic gradient descent (Adam) on the discretized cost $J^{N,\Delta t}$ of \eqref{eqn:J-cost-1} following \cite[Algorithm~20.1]{carmona2023deep}. In particular, we restrict the controls to a parametric neural class, pass to an $N$-particle system, discretize time, and run stochastic gradient descent on sampled costs, with the graphon control as the feedback policy. We note that this method requires direct access to the interaction function $\phi$. 

\textcolor{black}{Row-stochasticity and non-negativity hold exactly by \eqref{eqn:nn-softmax-labelonly}, and bounded scores $|s_{ij}|\le S$ give $N W_{\theta,ij}\le e^{2S}$. If the scores are also Lipschitz in the receiver label, with a common bound for all times, sender labels and discretization levels, the associated histogram kernels have receiver variation bounded independently of $N$; this bound depends on the network parameters. We do not impose a prescribed receiver-variation bound $\Vxi\le R$ during training, but instead measure the discrete variation \eqref{eqn:receiver-tv-discrete} of every learned kernel.}

Each experiment is run alongside the uncontrolled mean-field baseline $w \equiv 1$. We report the trained value $\widehat{J}_N$, the terminal variance $\Var(\vec{x}_N(T))$, the terminal mean $\overline{x}_N(T)$, and the peak measured receiver variation $\Vxi(\widehat{w}_N(t))$. Model-specific diagnostics include the cluster count for the bounded confidence model of \Cref{subsec:exp-bc}, in addition to the synchronization order parameter and winding number for the Kuramoto model of \Cref{subsec:exp-kuramoto}. Each configuration is run at a fixed and archived seed. The noise study of \Cref{subsec:exp-kuramoto} averages over $32$ independent Brownian paths. Configurations and code can be found in the repository \href{https://github.com/Hydrodynamical/Graphon_Control}{Graphon\_Control}.

\subsection{Cluster formation under bounded confidence}\label{subsec:exp-bc}

In this subsection, we study the bounded confidence model, where the interaction is a smoothed bounded-confidence influence. For a confidence radius $\varepsilon > 0$ and gain $a > 0$, define
\begin{equation}\label{eqn:phi-bc}
\phi_\varepsilon(z) := a \, z \left( 1 - \frac{z^2}{\varepsilon^2} \right)_+^{2}, \qquad (\,\cdot\,)_+ := \max(\,\cdot\,, 0).
\end{equation}
This is a $C^{1,1}$ mollification of the sharp Hegselmann--Krause influence $a\,z\,\mathbf{1}\{|z| \le \varepsilon\}$ \cite{Hegselmann-Krause}. 

We call a state \emph{fragmented} if its opinions split into internally consensual clusters with pairwise gaps exceeding $\varepsilon$. Because $\phi_\varepsilon$ vanishes outside $[-\varepsilon, \varepsilon]$, every pairwise force in a fragmented state is zero, so fragmented states are stationary points of the model dynamics \emph{for every admissible kernel}. Once the population polarizes, no admissible control can reconnect it. Therefore, a successful control must route influence through intermediate opinions while chains of overlapping confidence intervals still exist. We quantitatively track the cluster count by defining
$$
K(t) := 1 + \#\big\{ i : x_{(i+1)}(t) - x_{(i)}(t) > \varepsilon \big\},
$$
where $x_{(1)} \le \dots \le x_{(N)}$ is the sorted state.

We now move on to the experimental setup. We set $\varepsilon=1$ and $a=4$, and initialize the datum as $x(0, \xi) = s (\xi - \tfrac12) + A \sin (2\pi m \xi)$ with $s = 6$, $A = \frac{1}{2}$, $m = 4$, so the opinions spread over six confidence radii and are monotone in the label up to a small seeded ripple. We choose target $\xtar = 0$ with horizon $T = 10$, $\beta=1$ and $\gamma =100$. We describe the diagnostic numbers in \Cref{BC-control-comparison} below.

\begin{table}[ht]
\centering
\setlength{\tabcolsep}{10pt}
\renewcommand{\arraystretch}{1.3}
\begin{tabular}{|l|c|c|c|c|c|}
\hline
\bf Control & $K(T)$ &  $\overline{x}_N(T)$ & $\widehat{J}$ & $\Var(\vec{x}_N(T))$ & Peak $V_\xi(\widehat{w}(t))$\\
\hline
Mean-field & $4$ & $0.000$ & $38.08$ & $2.800$ &  0\\ 
\hline
Trained & $1$ & $-0.001$ & $17.93$ & $< 10^{-12}$ & 46.3\\
\hline
\end{tabular}
\caption{Terminal statistics for two types of controls in bounded confidence case.}
\label{BC-control-comparison}
\end{table}

\begin{figure}[ht]
    \centering
    \includegraphics[width=0.85\linewidth]{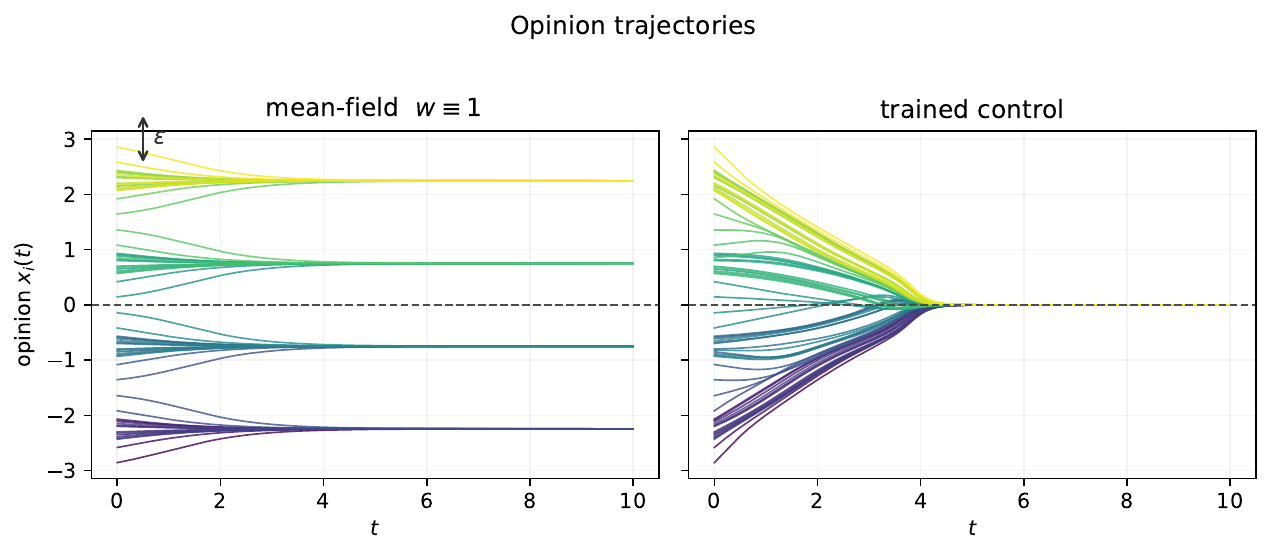}
    \caption{Opinion trajectories $x_i(t)$ under the uncontrolled mean-field dynamics and the trained control, colored by label $\xi$; the vertical bar shows the confidence radius $\varepsilon$.}
    \label{fig:bc-traj}
\end{figure}

\begin{figure}[ht]
    \centering
    \includegraphics[width=\linewidth]{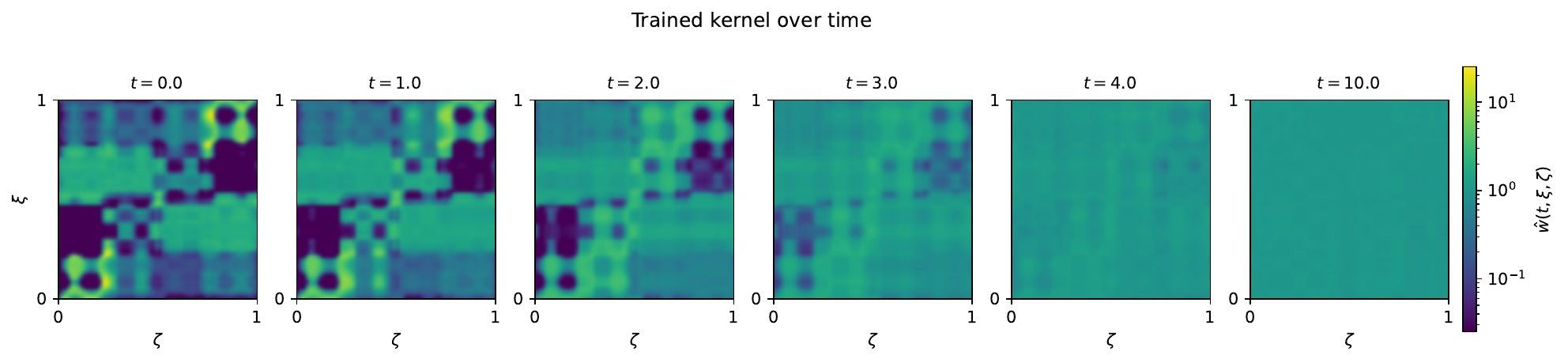}
    \caption{The trained kernel $\widehat{w}(t,\xi,\zeta)$ during the merge phase and at the horizon (log color scale): transient block-structured attention aligned with the nascent opinion clusters, relaxing to uniform after consensus.}
    \label{fig:bc-kernel}
\end{figure}

\Cref{fig:bc-traj} display the controlled and uncontrolled dynamics. In the uncontrolled case, the population fragments by $t \approx 1.9$, it has frozen into four equally sized clusters where every inter-cluster gap exceeds $\varepsilon$. The trained control instead funnels the population to consensus at the target by $t \approx 4$. \Cref{fig:bc-kernel} shows the \emph{transient} structure of the learned kernel. Once consensus is reached, the kernel relaxes toward uniform attention due to the regularizing term $\beta \Vert w \Vert^2_{L^2(I^2)}$.

\subsection{Untwisting twisted Kuramoto states}\label{subsec:exp-kuramoto}

In this subsection, we study the Kuramoto model where the interaction takes the form $\phi(z) = \sin z$. For identical natural frequencies, \eqref{eqn:graphon-dynamics} with $\phi = \sin$ \emph{is} the graphon Kuramoto model \cite{Kuramoto1975}.

Kuramoto introduced a notion of coherence by defining an order parameter \cite{Kuramoto1984}. For the continuum state $x(t, \cdot)$ and its $N$-agent discretization $\vec{x}_N(t)$, define
\begin{equation}\label{eqn:order-parameter}
r(t) := \Big| \int_I e^{i x(t,\xi)} \, d\xi \Big|,
\qquad
r_N(t) := \Big| \frac{1}{N} \sum_{j=1}^{N} e^{i x_N^j(t)} \Big|,
\end{equation}
the modulus of the mean of the state mapped to the unit circle, where $0 \le r \le 1$. This measures coherence as $r(t) = 1$ if and only if $x(t, \cdot)$ is almost everywhere constant modulo $2\pi$, and $r = 0$ on the $k$-twisted state $e^{i x(\xi)} = e^{2\pi i k (\xi - 1/2)}$ for any integer $k \ne 0$. We call these $k$-twisted states \emph{maximally incoherent}. 

To further investigate the behavior of states which differ by integer multiples of $2\pi$, we additionally introduce the winding number
\begin{equation}\label{eqn:winding}
\mathcal{W}(t) := \frac{1}{2\pi} \sum_{i = 1}^{N} \big[ x(t, \xi_{i+1}) - x(t, \xi_i) \big]_{2\pi},
\qquad
[d]_{2\pi} := \big( (d + \pi) \bmod 2\pi \big) - \pi,
\end{equation}
summing cyclically over the label grid ($\xi_{N+1} := \xi_1$). The winding number equals $k$ on the $k$-twisted state and $0$ at consensus.

We remark that twisted states are \emph{stationary points} of the uncontrolled dynamics for every symmetric convolution kernel. If $w(\xi, \zeta) = \rho(\xi - \zeta)$ with $\rho$ even and $x(\xi) = 2\pi k (\xi - \tfrac12)$ odd, then
$$
\int_I \rho(\xi - \zeta) \, \sin\big( 2\pi k (\zeta - \xi) \big) \, d\zeta = 0.
$$
This includes the mean-field baseline $\rho \equiv 1$, so the uncontrolled system sits at the twist forever without variance decay; for local ring couplings, twisted states are moreover \emph{stable} equilibria \cite{WileyStrogatzGirvan2006}. Therefore a successful controller has to break the translational symmetry that pins the equilibrium.

We now move on to the experimental setup. We initialize the system to be the $k$-twisted state $x(0, \xi) = 2\pi k (\xi - \tfrac12)$ with target $\xtar = 0$, and set $\beta =10$ and $\gamma = 100$.  The main experimental run take $k = 3$ at horizon $T = 20$ with $3000$ gradient steps, we summarize the diagnostic numbers in \Cref{Kuramoto-control-comparison}. 

\begin{table}[ht]
\centering
\setlength{\tabcolsep}{10pt}
\renewcommand{\arraystretch}{1.3}
\begin{tabular}{|l|c|c|c|c|}
\hline
\bf Control & $r(T)$ & $\widehat{J}$ & $\Var(\vec{x}_N(T))$ & Peak $V_\xi(\widehat{w}(t))$\\
\hline
Mean-field & $0.000$ & $792.03$ & $29.602$ & 0\\
\hline
Trained & $1.000$ & $368.62$ & $3.6 \times 10^{-8}$ & 24.4 \\
\hline
\end{tabular}
\caption{Terminal statistics for two types of controls in bounded Kuramoto case.}
\label{Kuramoto-control-comparison}
\end{table}

\begin{figure}[ht]
    \centering
    \includegraphics[width=0.85\linewidth]{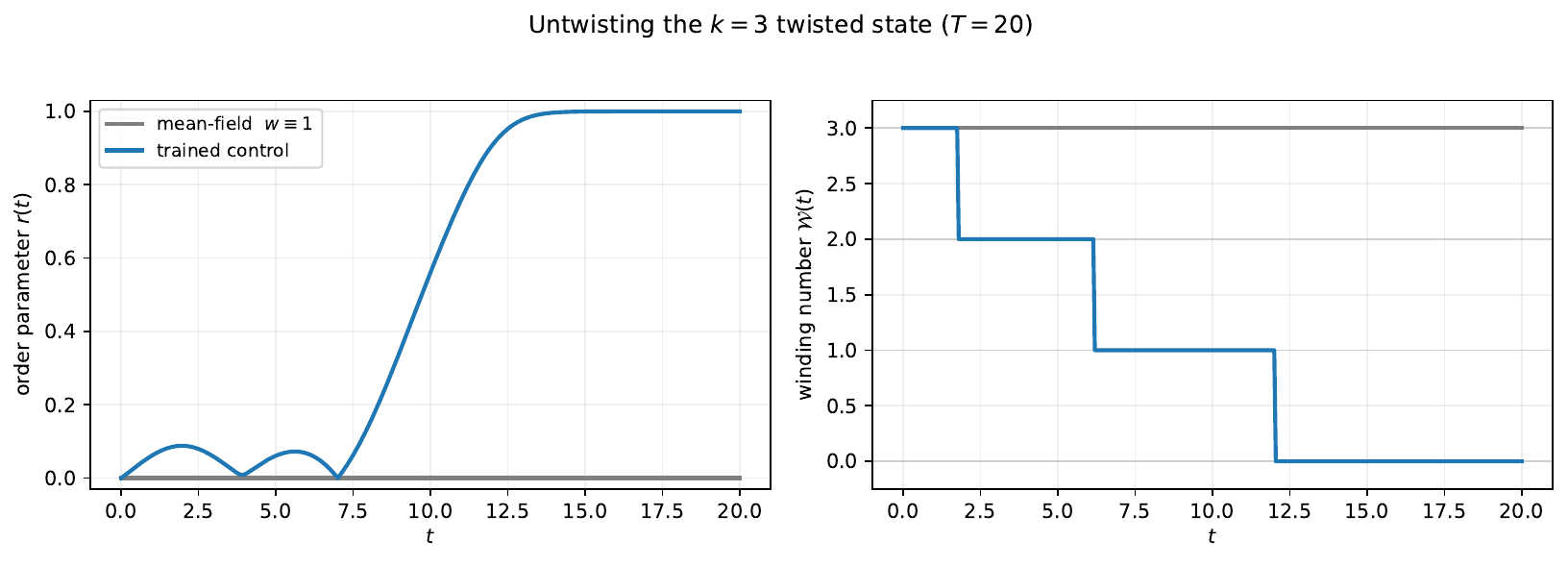}
    \caption{Order parameter $r(t)$ and winding number $\mathcal{W}(t)$ for the $k=3$ twisted state under the mean-field baseline and the trained control ($N = 64$).}
    \label{fig:kur-r}
\end{figure}

\Cref{fig:kur-r} displays the synchronization order parameter and winding number of the $k = 3$ run. The mean-field baseline remains maximally incoherent with $r(t) = 0$ for all time, with the variance frozen at its initial value $\pi^2 k^2/3 = 3\pi^2 \approx 29.6$. 

\begin{figure}[ht]
    \centering
    \includegraphics[width=\linewidth]{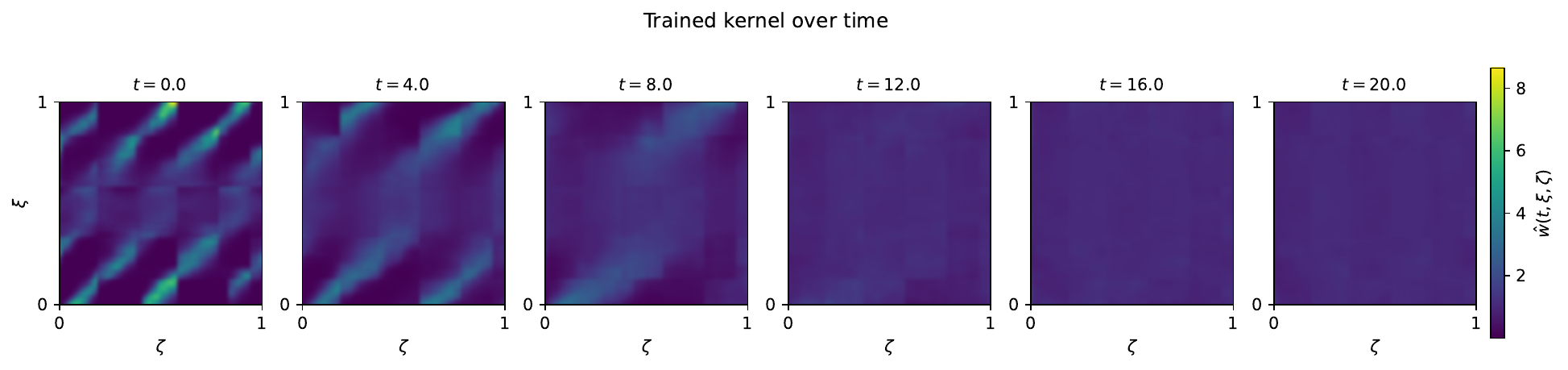}
    \caption{The trained kernel $\widehat{w}(t,\xi,\zeta)$ at six times ($k=3$): the tilted band breaks the translational symmetry pinning the twist, then relaxes to uniform attention after consensus.}
    \label{fig:kur-kernel}
\end{figure}

\Cref{fig:kur-kernel} shows the transient structure of the learned kernel. The trained control unwinds the twist completely. During the untwisting, attention concentrates along a band tilted off the diagonal, each receiver attends to senders at a fixed label offset. Once consensus is reached, the kernel relaxes to uniform due to the regularizing term $\beta \Vert w \Vert_{L^2}^2$.

Next, we compare twists $k = 1 (T = 5)$ and $k = 2 (T = 12)$, summarizing our results in \Cref{Different-k-comparison}. The $k=2$ problem is revisited under noise $\sigma > 0$ at the end of this subsection. We note the peak of $\Vxi(w(t))$ grows with depth, as deeper twists demand more receiver-heterogeneous kernels. 

\begin{table}[ht]
\centering
\setlength{\tabcolsep}{10pt}
\renewcommand{\arraystretch}{1.3}
\begin{tabular}{|l|c|c|c|c|c|}
\hline
$k$ & $T$ & $\widehat{J}_{\mathrm{trained}}$ & $\widehat{J}_{\mathrm{mean\text{-}field}}$ & Peak $\Vxi(\hat{w}(t))$ & $\Var(\vec{x}_N(T))$ \\
\hline
$1$ & $5$  & $57.52$  & $66.45$  & $3.7$  & $1.0 \times 10^{-2}$ \\
\hline
$2$ & $12$ & $176.09$ & $277.88$ & $12.0$ & $1.5 \times 10^{-5}$ \\
\hline
$3$ & $20$ & $368.62$ & $792.03$ & $24.4$ & $3.6 \times 10^{-8}$ \\
\hline
\end{tabular}
\caption{Trained and mean-field costs for Kuramoto twists of increasing depth, with terminal variance and peak receiver variation of the trained control.}
\label{Different-k-comparison}
\end{table}

\begin{figure}[ht]
    \centering
    \includegraphics[width=0.85\linewidth]{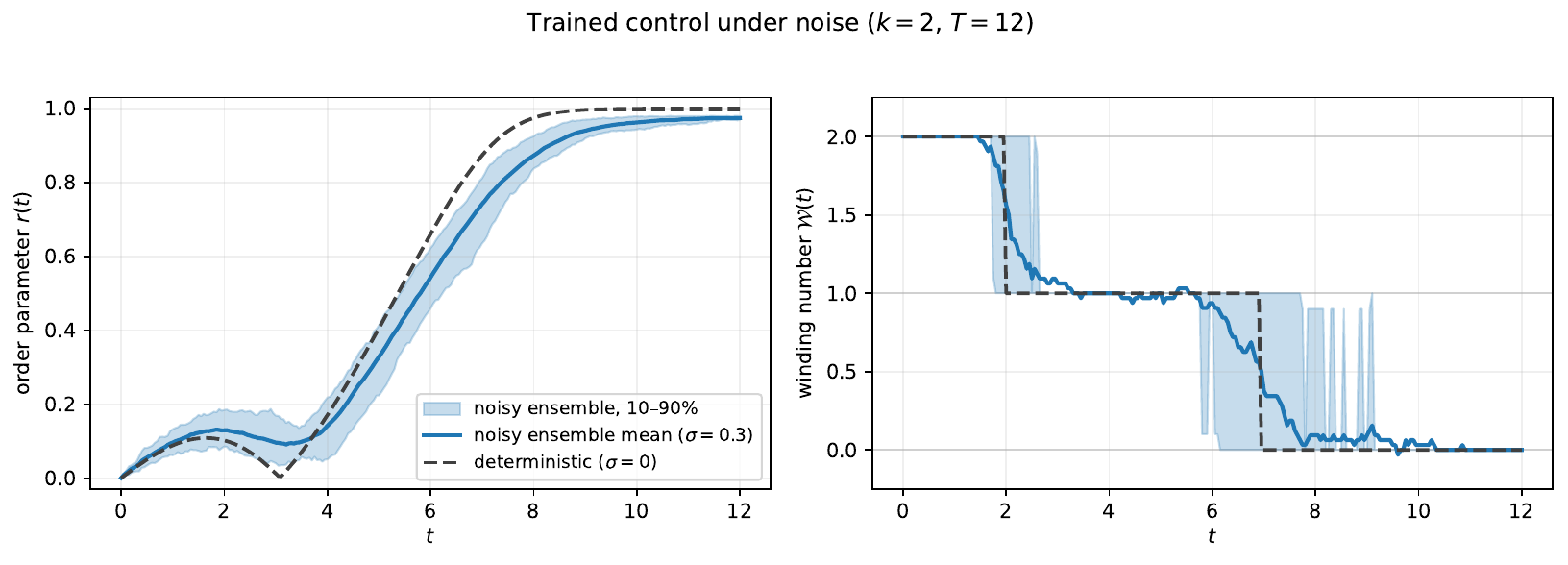}
    \caption{The $k = 2$ control problem under Brownian noise ($\sigma = 0.3$): order parameter (left) and winding number (right) for the control trained on the noisy dynamics, evaluated over $32$ Brownian paths (mean and $10$--$90\%$ band) against the deterministic run.}
    \label{fig:kur-noise}
\end{figure}

While the theory of \Cref{sec:existence} is deterministic, the method of \Cref{sec:numerical-method} allows for random dynamics. To investigate robustness against noise, we retrain the $k = 2$, $T = 12$ problem with Brownian noise $\sigma = 0.3$ in the training dynamics, and evaluate the controlled dynamics over $32$ Brownian paths started at the exact twist. \Cref{fig:kur-noise} compares the ensemble statistics of $r(t)$ and $\mathcal{W}(t)$ with the deterministic run. Every path ends at $\mathcal{W} = 0$, while the terminal order parameter concentrated at mean $r(T) = 0.974$.

We deliberately limited the numerics to these two experiments. The framework invites others, which we leave to future work, such as alignment models in which classical distance-based communication rates become the baseline to beat, \textcolor{black}{the selection effect of the regularizer on the example of \Cref{rem:non-uniqueness}}, a quantitative study of the discrete-to-continuum convergence asserted by \Cref{thm:main}, and a systematic treatment of the stochastic dynamics beyond the robustness check of \Cref{subsec:exp-kuramoto}.

\end{document}